\documentclass[]{interact}

\usepackage{epstopdf}
\usepackage[caption=false]{subfig}

\usepackage[numbers,sort&compress]{natbib}
\bibpunct[, ]{[}{]}{,}{n}{,}{,}
\renewcommand\bibfont{\fontsize{10}{12}\selectfont}

\theoremstyle{plain}
\newtheorem{theorem}{Theorem}[section]
\newtheorem{lemma}[theorem]{Lemma}

\theoremstyle{definition}

\theoremstyle{remark}

\usepackage{algorithm,algpseudocode}
\usepackage{hyperref}
\usepackage{multirow}

\newcommand{\R}{\mathbb{R}}

\newcommand{\E}{\mathbb{E}}

\DeclareMathOperator{\Var}{Var}
\DeclareMathOperator{\trace}{trace}
\DeclareMathOperator{\Hess}{Hess}
\DeclareMathOperator{\sign}{sign}

\begin{document}


\title{A Smoothing Algorithm for $\ell^{1}$ Support Vector Machines\thanks{This work is an extension of \emph{Smoothed Hinge Loss and $\ell^{1}$ Support Vector Machines}, published in International Conference on
Data Mining Workshops, held in Singapore, November 17-20, 2018.}}

\author{
\name{Ibrahim Emirahmetoglu\textsuperscript{a}, Jeffrey Hajewski\textsuperscript{b}, Suely Oliveira\textsuperscript{b}, and David E.~Stewart\textsuperscript{a}}
\affil{\textsuperscript{a}Department of Mathematics, The University of Iowa; \textsuperscript{b}Department of Computer Science, The University of Iowa}
}

\maketitle

\begin{abstract}
A smoothing algorithm is presented for solving the soft-margin Support Vector Machine (SVM) optimization problem with an $\ell^{1}$ penalty. This algorithm is designed to require a modest number of passes over the data, which is an important measure of its cost for very large datasets. The algorithm uses smoothing for the hinge-loss function, and an active set approach for the $\ell^{1}$ penalty. The smoothing parameter $\alpha$ is initially large, but typically halved when the smoothed problem is solved to sufficient accuracy. Convergence theory is presented that shows $\mathcal{O}(1+\log(1+\log_+(1/\alpha)))$ guarded Newton steps for each value of $\alpha$ except for asymptotic bands $\alpha=\Theta(1)$ and $\alpha=\Theta(1/N)$, with only one Newton step provided $\eta\alpha\gg1/N$, where $N$ is the number of data points and the stopping criterion that the predicted reduction is less than $\eta\alpha$. The experimental results show that our algorithm is capable of strong test accuracy without sacrificing training speed.
\end{abstract}

\begin{amscode}
Primary: 65K05, Secondary: 49M15
\end{amscode}

\section{Introduction} \label{sec:Intro}
Dealing with large datasets has lead to a strong interest in methods that have low iteration costs, such as stochastic gradient descent (SGD) \cite{Niu2011, Rob51}. More classical methods such as Newton's method for optimization \cite[\S3.3]{Noc2006} are generally not used as their cost per iteration involves solving linear systems which takes $\mathcal{O}(m^{3}+Nm^{2})$ operations where $m$ is the number of unknowns and $N$ the number of data items. Contrary to conventional wisdom, we argue that Newton's method and more sophisticated line search methods are actually more appropriate for very large data problems, since their computational issues are typically due to the large number of data items ($N$) rather than the dimension of the problem ($m$). Wide data, where the dimension $m$ of the data vectors $\bm{x}_{i}$ is large compared to the number of data items $N$, is still problematic for Newton's method as the Hessian matrix is then singular. However, in this paper we focus on $\ell^{1}$ Support Vector Machines ($\ell^{1}$\,SVMs) and argue that Newton's method with a suitable line search and an active-set strategy can also solve these problems very efficiently. The algorithm developed here is, in part, inspired by \cite{Osb2000} for the basis pursuit noise-reduction problem. A different version of this algorithm is reported in the conference paper \cite{Haj2018}. The development of the line search algorithm, the convergence analysis of Section~\ref{sec:convergence-rate} are not reported in the conference paper.

In what follows, the norm used is the 2-norm $\left\| \bm{z}\right\|_{2}=\sqrt{\bm{z}^{T}\bm{z}}$ unless otherwise indicated. The soft-margin SVM \cite[p.~263]{Dud2000} for given data $(\bm{x}_{i},y_{i})$, $i=1,\,2,\,\ldots,\,N$ where $y_{i}=\pm1$ for all~$i$, minimizes
\begin{equation} \label{eq:SVM-soft-margin}
\widehat{f}(\bm{w}; \mathcal{D}) := \frac{1}{2} \lambda\left\| \bm{w} \right\|^{2} + \frac{1}{N}\sum_{i=1}^{N} \max(0,\,1-y_{i}\,\bm{w}^{T}\bm{x}_{i})
\end{equation}
over all $\bm{w} \in \mathbb{R}^{m}$. Here the value of $\lambda\geq0$ is used to control the norm of the vector $\bm{w}$ when $\lambda>0$. The function $\max(0,\,1-y_{i}\,\bm{w}^{T}\bm{x}_{i})$ is called the \emph{hinge-loss function} as it is based on the function $u\mapsto\max(0,\,u)$ whose graph looks like a hinge. The $\ell^{1}$\,SVM
for the same data minimizes
\begin{equation} \label{eq:L1-SVM}
f(\bm{w}; \mathcal{D}) := \frac{1}{2}\lambda\left\Vert \bm{w}\right\Vert ^{2}+\frac{1}{N}\sum_{i=1}^{N}\max(0,\,1-y_{i}\,\bm{w}^{T}\bm{x}_{i})+\mu\left\Vert \bm{w}\right\Vert _{1}
\end{equation}
over $\bm{w}$. Here $\mu>0$ controls the level of sparsity of $\bm{w}$. Larger values tend to mean fewer components of $\bm{w}$ are non-zero; if $\mu$ is large enough then $\bm{w}=0$. Note that this formulation is similar to, but not the same as the 1-norm SVM of Zhu, Rosset, Hastie and Tibshirani \cite{Zhu2003}. Also, the algorithm obtained here takes $\mathcal{O}(N)$ time per step with respect to the number of data points, while the algorithm of Zhu et~al.\ is $\Omega(N^{2})$ as it involves identifying the intersections of a descent line with the hyperplanes $1-y_{i}\,\bm{w}^{T}\bm{x}_{i}=0$ for each data point. Rather we use a smoothing approach for the sum of the hinge-loss functions. 

Traditionally, for optimization problems, the numbers of function, gradient, and Hessian matrix evaluations are used to measure the cost of the algorithm. For large-scale data mining types of optimization problems, perhaps a different measure of performance is more important: the number of passes over the data. The general form of most optimization problems used in data mining is
\begin{equation} \label{eq:data-mining-gen-form}
\min_{\bm{w}}f(\bm{w}; \mathcal{D}):=R(\bm{w})+\frac{1}{N}\sum_{i=1}^{N}\varphi(\bm{x}_{i},y_{i};\bm{w})
\end{equation}
where $\mathcal{D}=\left\{ \,(\bm{x}_{i},y_{i})\mid i=1,\ldots,N\,\right\} $ is the dataset, $R$ is a regularization function, and $\varphi(\bm{x}, y; \bm{w}) := \max(0,\,1-y\,\bm{w}^{T}\bm{x})$ is the hinge loss function. Provided $\bm{w}$ has relatively low dimension (say, below $10^{3}$) and $N$ is large (say, $10^{5}$ to $10^{9}$), the cost of computing $R(\bm{w})$ is modest and can be computed on one processor, while the computations of $\varphi(\bm{x}_{i}, y_{i}; \bm{w})$ should be carried out in parallel, and then summed via a parallel reduction operation \cite{Ble96}.

Computing the gradient of the objective function
\[
\nabla f(\bm{w}; \mathcal{D}) = \nabla R(\bm{w}) + \frac{1}{N} \sum_{i=1}^{N} \nabla_{\bm{w}} \varphi(\bm{x}_{i},y_{i};\bm{w}),
\]
can be carried out in a similar manner to the objective function, except that the reduction (summation) is applied to the gradients $\nabla_{\bm{w}} \varphi(\bm{x}_{i}, y_{i}; \bm{w})$. Similarly, the Hessian matrices can be computed in parallel but the reduction (summation) is applied to the Hessian matrices $\Hess_{\bm{w}} \varphi(\bm{x}_{i},y_{i};\bm{w})$ of the loss functions. This may become an expensive step if $m$ becomes large, as the reduction must be applied to objects of size $\mathcal{O}(m^{2})$ where $\bm{w} \in \R^{m}$. In such cases, a BFGS algorithm may be appropriate instead of a direct Newton method. For the $\ell^{1}$\,SVM problem (\ref{eq:L1-SVM}), however, for modest $\mu>0$ there are relatively few components of $\bm{w}$ that are non-zero at the optimum.

The use of the $\ell^{1}$ penalty term $\mu\left\| \bm{w}\right\|_{1}$ with suitable $\mu>0$ results in sparse solutions; that is, the number of non-zero components of $\bm{w}$ can be small compared to $m$. This can be beneficial for a number of reasons, both statistical and computational. Here we will focus on the computational advantages. If $\left|\left\{ \,j\mid w_{j}\neq0\,\right\} \right|$ is small compared to $m$, then solving the linear system in the non-zero components of $\bm{w}$ can be much less computationally intensive than solving the system for all components of $\bm{w}$. However, this approach requires identifying the non-zero components of $\bm{w}$, and explicitly changing this ``active set'' when appropriate. This has consequences for other aspects of the algorithm developed here: the $\ell^{1}$ penalty should \emph{not} be smoothed. Instead we keep it in its original form, and use the non-smoothness to keep the active set small. In turn, this impacts other aspects of the algorithm, such as the line search method.

The method proposed in this paper must deal with several inter-related issues:
\begin{itemize}
\item The $\ell^{1}$ penalty enforces sparsity of solutions, but introduces non-smoothness into the objective function that is not smoothed in the algorithm. This is handled by using an ``active set'' method \cite[pp.~265--275]{Fle87}, and modifying the line search strategy.

\item The unregularized objective function $N^{-1} \sum_{i=1}^{N} \varphi(\bm{x}_{i},y_{i};\bm{w})$ is actually close to being smooth for large $N$, but direct computations of Hessian matrices of this function give useless results (either zero or undefined). Smoothing is used to give useful Hessian estimates.

\item The smoothing parameter should go to zero as the computation proceeds to accurately solve the intended problem (\ref{eq:L1-SVM}). Each reduction of the smoothing parameter should require at most $\mathcal{O}(1)$ Newton steps. We are able to show $\mathcal{O}(1+\log\log(1/\alpha))$ guarded Newton steps for $1\gg\alpha$, with most steps in the asymptotic range $1\gg\alpha\gg N^{-1}$ taking exactly one Newton step and no line search.
\end{itemize}
The remainder of the paper is as follows: Section~\ref{sec:Algorithm-development} shows how the algorithm is developed; Section~\ref{sec:convergence-rate} develops the theory arguing for a bounded number of Newton steps for each reduction of the smoothing parameter; Section~\ref{sec:results} gives results showing the practicality and competitiveness of the algorithm developed here. Finally, Section \ref{sec:Conclusion} concludes.

\section{Algorithm development} \label{sec:Algorithm-development}
In this section, we first describe smoothing the hinge-loss function and Hessian matrices. Next, we explain how the algorithm is developed and provide the pseudo-code for our smoothing algorithm. Finally, the issues with the line search procedure are addressed.

\subsection{Smoothing the hinge-loss function and Hessian matrices} \label{subsec:Smoothed-hinge-loss}
The hinge-loss function is defined by
\begin{equation} \label{eq:hinge-loss}
\varphi(\bm{x},y;\bm{w}) := \max(0,\,1-y\,\bm{w}^{T}\bm{x}),
\end{equation}
where $y \in \{\pm 1\}$ and $\bm{x}, \bm{w} \in \R^m$. Note that the hinge loss is in fact of the form 
\begin{equation} \label{eq:hinge-loss-2}
\varphi(\bm{x},y;\bm{w},b) = \max\big(0,\,1-y(\bm{w}^{T}\bm{x} + b)\big),
\end{equation}
where $b \in \R$ represents the bias term. However, defining $\widetilde{\bm{x}} = (\bm{x}, 1)$ and $\widetilde{\bm{w}} = (\bm{w}, b)$ in $\R^{m+1}$, the expression \eqref{eq:hinge-loss-2} can be written of the form \eqref{eq:hinge-loss}. For the sake of brevity, we stick with the form in \eqref{eq:hinge-loss}. 

\begin{figure}
\centering
\includegraphics[scale=0.55]{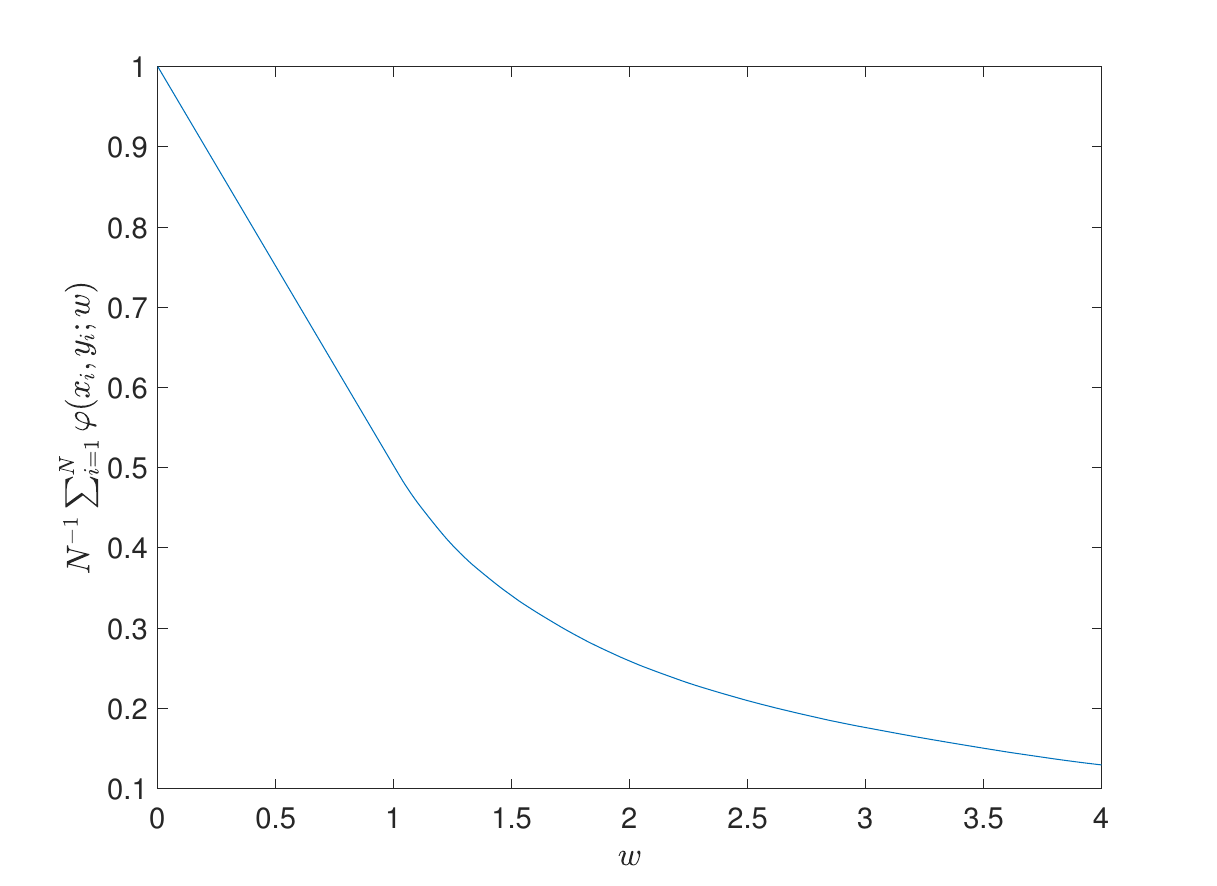}
\caption{\label{fig:nearly-smooth}Plot of $N^{-1}\sum_{i=1}^{N}\varphi(x_{i},y_{i};w)$
for randomly chosen data ($n=200$)}
\end{figure}

The hinge-loss function $\varphi(\bm{x},y;\bm{w})$ is a piecewise linear function of $\bm{w}$, and so its Hessian matrix is zero or undefined. Thus, $\frac{1}{N}\sum_{i=1}^{N}\varphi(\bm{x}_{i},y_{i};\bm{w})$ is also a piecewise linear function of $\bm{w}$, and hence its Hessian matrix is either zero or undefined. On the other hand, $\frac{1}{N}\sum_{i=1}^{N}\varphi(\bm{x}_{i},y_{i};\bm{w})$ usually appears to be very smooth. For example, if $N=200$, and for $y=+1$, $x$ chosen randomly and uniformly from $[0,+1]$ while for $y=-1$, $x$  chosen uniformly and randomly from $[-1,0]$, the function $\frac{1}{N}\sum_{i=1}^{N}\varphi(x_{i},y_{i};w)$ looks like Figure~\ref{fig:nearly-smooth}.

Under some statistical assumptions detailed below, the function $\frac{1}{N}\sum_{i=1}^{N} \varphi(\bm{x}_{i},y_{i};\bm{w})$ approaches a smooth function $h(\bm{w})$ as $N\to\infty$. Rather than computing the exact Hessian matrix of $\frac{1}{N}\sum_{i=1}^{N}\varphi(\bm{x}_{i},y_{i};\bm{w})$ with respect to $\bm{w}$ (which is either zero or undefined), we compute an approximation to the Hessian matrix of $h(\bm{w})$. This can be done by means of a smoothed hinge-loss function. We can approximate $\varphi(\bm{x},y;\bm{w})$ by a smoothed hinge-loss function $\varphi_{\alpha}(\bm{x},y;\bm{w})$ given by
\begin{equation} \label{eq:def-smoothed-hinge-loss}
\varphi_{\alpha}(\bm{x},y;\bm{w}) := \alpha\,\psi(u/\alpha) = \frac{1}{2}\left(u+\sqrt{\alpha^{2}+u^{2}} \right),
\end{equation}
where $u = 1-y\,\bm{w}^{T}\bm{x}$ and $\psi(u) := \frac{1}{2}\big( u + \sqrt{1 + u^2} \big)$. Using a smoothed hinge-loss function does not change the value of $N^{-1}\sum_{i=1}^{N}\varphi(\bm{x}_{i},y_{i};\bm{w})$ significantly, but does enable us to estimate the Hessian matrix of $h(\bm{w})$, as well as its gradient.

For large datasets which come from some statistical distribution with a $C^{1}$ probability density function, the mean of the hinge-loss functions approaches a $C^{2}$ function
\begin{equation} \label{eq:def-h(w)}
\frac{1}{N} \sum_{i=1}^{N} \varphi(\bm{x}_{i},y_{i};\bm{w}) \to \sum_{y\in\left\{ \pm1\right\} } \int_{\R^{m}}\varphi(\bm{x},y;\bm{w})\,p(\bm{x},y)\,d\bm{x} =: h(\bm{w})
\end{equation}
as $N\to\infty$. Here $p(\bm{x},y)$ is the probability density function of $(\bm{x},y)$. The difficulty is in estimating the Hessian matrix of this unknown smooth function. What we want is that 
\begin{equation} \label{eq:approx-Hess-h(w)}
\frac{1}{N}\sum_{i=1}^{N}\Hess_{\bm{w}}\varphi_{\alpha}(\bm{x}_{i},y_{i};\bm{w}) \approx \Hess\,h(\bm{w})
\end{equation}
for $N$ sufficiently large. Since $\Hess_{\bm{w}}\varphi(\bm{x},y;\bm{w}) = \bm{x}\bm{x}^{T}\,\delta(1-y\,\bm{w}^{T}\bm{x})$ where $\delta$ is the Dirac-delta distribution, we can write
\begin{align*}
\Hess h(\bm{w}) &= \sum_{y\in\left\{ \pm1\right\} }\int_{\R^{m}} \Hess_{\bm{w}}\varphi(\bm{x},y; \bm{w})\,p(\bm{x},y)\,d\bm{x} \\
&= \frac{1}{\|\bm{w}\|} \sum_{y\in\left\{ \pm1\right\} }\int_{\left\{ \bm{x}\mid y=\bm{w}^{T}\bm{x}\right\} }\bm{x}\bm{x}^{T}p(\bm{x},y)\,dS(\bm{x}).
\end{align*}
On the other hand, we have $\Hess_{\bm{w}} \varphi_\alpha(\bm{x},y; \bm{w}) \to \bm{x}\bm{x}^T \delta(1-y\,\bm{w}^T\bm{x})$ as $\alpha \downarrow 0$, since
\begin{align*}
\frac{d^{2}}{du^{2}}\alpha\psi(u/\alpha) & =\frac{1}{\alpha}\psi''(u/\alpha)=\frac{1}{2}\frac{\alpha^{2}}{(\alpha^{2}+u^{2})^{3/2}}\to\delta(u)\quad \text{ as }\alpha \downarrow 0,
\end{align*}
in the sense of distributions and weakly{*} in the sense of measures. Thus, for continuous $p$, we have
\begin{align*}
\sum_{y\in\left\{ \pm1\right\} }\int_{\R^{m}}\Hess_{\bm{w}}\varphi_{\alpha}(\bm{x},y;\bm{w})\,p(\bm{x},y)\,d\bm{x} 
&\to \frac{1}{\|\bm{w}\|} \sum_{y\in\left\{ \pm1\right\} }\int_{\left\{ \bm{x}\mid y=\bm{w}^{T}\bm{x}\right\} }\bm{x}\bm{x}^{T}p(\bm{x},y)\,dS(\bm{x}) \\
&= \Hess h(\bm{w})
\end{align*}
as $\alpha \downarrow 0$. Moreover, if we choose $(\bm{x}_{i},y_{i})$ independently, distributed according to the probability distribution $p$, and the variance for the probability distribution $p$ is finite, then by the Strong Law of Large Numbers \cite[p.~239]{Loe77},
\[
\frac{1}{N}\sum_{i=1}^{N}\Hess_{\bm{w}}\varphi_{\alpha}(\bm{x}_{i},y_{i};\bm{w})\to\sum_{y\in\left\{ \pm1\right\} }\int_{\R^{m}}\Hess_{\bm{w}}\varphi_{\alpha}(\bm{x},y;\bm{w})\,p(\bm{x},y)\,d\bm{x}
\]
as $N \to \infty$ almost surely. To make this work in a practical sense, we need the number of samples $N$ to be ``sufficiently large'' for a given $\alpha>0$ in order to have \eqref{eq:approx-Hess-h(w)} at least with high probability. A natural question is how large $N$ has to be for a given $\alpha$ in order to have a good approximation. Since $\varphi_{\alpha}(\bm{x},y;\bm{w})$ only depends on $(1-y\bm{w}^{T}\bm{x})/\alpha$, we only need $N\alpha\gg\left\Vert \bm{w}\right\Vert $ as $\alpha\downarrow0$ in order to achieve a given level of accuracy in approximating $\Hess h(\bm{w})$. Hence, the number of data points needed to obtain a good approximation of the curvature of the objective function is not exorbitant. This is made more rigorous in Section~\ref{sec:convergence-rate}.

As noted above, $\varphi(\bm{x}_{i},y_{i};\bm{w})$ is piecewise linear in $\bm{w}$, so $\Hess f(\bm{w}; \mathcal{D})$ is either undefined or $\Hess R(\bm{w})$ where this Hessian is defined. Using this for a Newton method will not give fast convergence. Instead, for large $N$ we wish to get a suitable approximation to $\Hess(h(\bm{w})+R(\bm{w}))$, while taking into account any non-smoothness in $R(\bm{w})$. The problem then is to choose $\alpha$. From the analysis above, we could choose $\alpha$ to be inversely proportional to $N$.

Instead, we use an approach which is agnostic regarding $N$ and the distribution of data points: start with a large value of $\alpha$, minimize $f_{\alpha}(\bm{w}; \mathcal{D})$ over $\bm{w}$, then repeatedly reduce $\alpha$ by (for example) halving $\alpha$, minimizing $f_{\alpha}(\bm{w}; \mathcal{D})$ over $\bm{w}$ for each $\alpha$. This can still be a very efficient algorithm, with $\mathcal{O}(1+\log(1+\log_+(1/\alpha)))$ Newton steps as $\alpha \downarrow 0$ typically being taken for each value of $\alpha>0$, often using just one Newton step for each value of $\alpha$, as discussed in Section~\ref{sec:convergence-rate}. Note that $\log_+(u) := \max(0,\,\log u)$. 

\subsection{Smoothing algorithm} \label{subsec:Smoothing-algorithm}
Replacing a hinge-loss function with a smoothed hinge-loss $\varphi_\alpha(\bm{x},y; \bm{w})$ given in \eqref{eq:def-smoothed-hinge-loss}, we define the following objective functions for $\ell^1$\,SVM and $\ell^1$--$\ell^2$\,SVM, respectively:
\begin{align*}
\widehat{f}_{\alpha}(\bm{w}; \mathcal{D}) &= \frac{1}{2}\lambda\bm{w}^{T}\bm{w} + \frac{1}{N} \sum_{i=1}^{N} \varphi_{\alpha}(\bm{x}_i, y_{i}; \bm{w}), \\
f_{\alpha}(\bm{w}; \mathcal{D}) & =\widehat{f}_{\alpha}(\bm{w}; \mathcal{D})+\mu\left\Vert \bm{w}\right\Vert_{1}.
\end{align*}
Note that $\nabla\widehat{f}_{\alpha}(\bm{w})$ is well-defined for all $\bm{w} \in \R^m$ provided $\alpha>0$, but that $f_{\alpha}$ is not smooth. A high-level algorithm is outlined in Algorithm~\ref{alg:SVM-l1-penalty-1-1}. 

\begin{algorithm}[H]
\begin{algorithmic}[1]
\Require{$\alpha,\,\alpha_{\min},\,\beta,\,\mu,\,\lambda>0$}
\Function{SVMsmooth}{$X,\,\bm{y},\,\bm{w},\,\lambda,\,\mu,\,\alpha,\,\alpha_{\min}, \beta$}
\State$\mathcal{I}\gets\left\{ \,j\mid w_{j}=0\,\right\}$; \quad $\overline{\mathcal{I}} \gets \left\{ \,j \mid w_j \neq 0 \right\}$ \Comment{Inactive\,/\,active sets}
\State$\mathcal{J}\gets\left\{ \,j\in\mathcal{I}\mid\left| \widehat{g}_{j}\right|>\mu\,\right\} $\Comment{Add to active set}
\State$\mathcal{I} \gets \mathcal{I} \setminus \mathcal{J}$; \quad $\overline{\mathcal{I}} \gets \overline{\mathcal{I}} \cup \mathcal{J}$
\State$adjust \gets \texttt{true}$
\While{$\alpha>\beta\,\alpha_{\min}$}\Comment{While smoothing parameter not at threshold}
  \State{}Carry out Newton steps on smoothed problem (Algorithm~\ref{alg:smoothed-l1-Newton-step})
\EndWhile
\State$\mathbf{return}\;\bm{w}$
\EndFunction
\end{algorithmic}
\caption{Smoothing algorithm for SVM with $\ell^{1}$ penalty}
\label{alg:SVM-l1-penalty-1-1}
\end{algorithm}
For the implementation of our algorithm, it is convenient to write the dataset ${\mathcal{D}=\left\{(\bm{x}_{i},y_{i})\mid i=1,\ldots,N\right\}}$ in matrix-vector form as $\mathcal{D} = (X, \bm{y})$, where
\begin{align*}
X := \left[ \begin{array}{cccc}
\bm{x}_1 & \bm{x}_2 & \cdots & \bm{x}_N
\end{array} \right]^T \in \R^{N \times m}, \qquad
\bm{y} := \left( y_1,\ y_2,\ \ldots,\ y_N  \right) \in \R^{N},
\end{align*}
so that the inputs $u_i := 1-y_{i}\,\bm{x}_{i}^{T}\bm{w}$ to $\psi(u_i/\alpha)$ form the vector $\bm{u} := \bm{e}-\bm{y}\circ(X\bm{w})$. Note that $\circ$ represents Hadamard product, and $\bm{e}$ is the vector of 1's in $\R^N$. This vector formulation is helpful in languages such as Matlab$^{TM}$ and Julia.

The algorithm used can be broken down into a number of phases. At the top level, the method can be considered as applying Newton's method to a smoothed problem (smoothing parameter $\alpha$) keeping an active set $\overline{\mathcal{I}} = \left\{ \,j\mid w_{j} \neq 0\,\right\} $. The active set will need to change, either by gaining elements where $w_{j} = 0$ but the gradient component $\widehat{g}_{j}= \frac{\partial\widehat{f}_{\alpha}}{\partial w_{j}}(\bm{w})$ satisfies $\left| \widehat{g}_{j}\right|>\mu$ indicating that allowing $w_{j}\neq0$ will result in a lower objective function value, or by losing elements where $w_{j}\neq0$ but $w_{j}+sd_{j}=0$ or $\sign(w_j) =-\sign(w_j+sd_j)$ resulting from the line search procedure. 

An essential choice in this algorithm is \emph{not} to smooth the $\ell^{1}$ penalty term, and instead use an active/inactive set approach. If we had chosen to smooth the $\ell^{1}$ penalty term, then the computational benefits of the smaller linear system in the Newton step $\bm{d}_{\overline{\mathcal{I}}}\gets-H_{\overline{\mathcal{I}},\overline{\mathcal{I}}}^{-1}{\bm{g}}_{\overline{\mathcal{I}}}$ would be lost. Instead, smoothing the $\ell^{1}$ term would mean that the linear system to be solved would have size $m\times m$, where $m$ is the dimension of $\bm{w}$. This would be particularly important for problems with wide datasets where $m$ can be very large. Instead, we expect that there would be bounds on the size of $\overline{\mathcal{I}}$, the number of active weights $w_{j}\neq 0$.

\begin{algorithm}[H]
\begin{algorithmic}[1]
  \State$\widehat{\bm{g}} \gets\nabla \widehat{f}_{\alpha}(\bm{w})$
  \State$\bm{g} \gets \widehat{\bm{g}} + \mu \sign(\bm{w})$
  \State$g_{j} \gets g_{j} - \mu \sign(\widehat{g}_{j})$ \quad $\forall\, j\in\mathcal{J}$
  \State$H \gets \Hess_{\bm{w}} \widehat{f}_{\alpha}(\bm{w})$
  \State$\bm{d}_{\overline{\mathcal{I}}}\gets-H_{\overline{\mathcal{I}},\overline{\mathcal{I}}}^{-1}\bm{g}_{\overline{\mathcal{I}}}$; \quad $\bm{d}_{\mathcal{I}}\gets \bm{0}$; \Comment{full Newton step}

  \If{$\big|\bm{d}^{T}\bm{g}\big| < \eta\,\alpha$}\Comment{If smoothed problem nearly solved for $\alpha$ and $\overline{\mathcal{I}}$}
    \State{}Adjust active set \& reduce smoothing parameter (Algorithm~\ref{alg:Adjust-active-set})
  \EndIf
  \State$s\gets\text{\sc{LinesearchL1}}(\bm{w},\bm{d},\bm{d}^{T}\widehat{\bm{g}},\frac{1}{2}\bm{d}^{T}H\bm{d},\mu,s_{\max})$
  \State $\bm{w}^{+}\gets\bm{w}+s\,\bm{d}$

  \While{$f_{\alpha}(\bm{w}^{+})>f_{\alpha}(\bm{w})+c_{1}\,s\,\bm{d}^{T}\bm{g}$} \Comment{Armijo line search}
    \State $\mathbf{if}$ $s < s_{\min}$, $\mathbf{then}$ $\mathbf{break}$ \Comment{Line search failure}
    \State$s\gets s/2$; \quad $\bm{w}^{+}\gets\bm{w}+s\,\bm{d}$;
  \EndWhile
  \If{$s < s_{\min}$} \Comment{If line search failed}
    \State $\alpha \gets \beta\, \alpha$ \Comment{Reduce $\alpha$ and optimize for this new $\alpha$}
    \State $adjust \gets \texttt{true}$
    \State $\mathbf{continue}$
  \EndIf
  \State $\mathcal{I} \gets \big\{ \,j \mid w_{j}^+=0 \ \text{or}\ \sign(w_j^+) = -\sign(w_j) \,\big\}$ \Comment{Add to $\mathcal{I}$ if line search indicates}
  \State${w}_j \gets \begin{cases} 0 & \text{if } j \in \mathcal{I}, \\[-0.3em] {w}_j^{+} & \text{otherwise.} \end{cases}$
\end{algorithmic}
\caption{Newton step}
\label{alg:smoothed-l1-Newton-step}
\end{algorithm}

The Newton step computations are shown in Algorithm~\ref{alg:smoothed-l1-Newton-step}. We first compute the gradient and the Hessian matrix. Care must be taken at this point to ensure that we compute the correct gradient for the components $j\in\mathcal{J}$ with $w_{j}=0$ but $\left|\widehat{g}_{j}\right|>\mu$. In this case, we need to add $j$ to the active set, and we will increase $w_{j}$ if $\widehat{g}_{j}<0$ and decrease $w_{j}$ if $\widehat{g}_{j}>0$. The full Hessian matrix is not actually needed, just the ``active'' part of the Hessian matrix: $H_{\overline{\mathcal{I}},\overline{\mathcal{I}}}$. The Newton step $\bm{d}$ is then computed. If the predicted reduction of the function value is sufficiently small, we can then assume the problem for the current active set $\overline{\mathcal{I}}$ and smoothing parameter $\alpha>0$ has been solved to sufficient accuracy. Next, we can either adjust the current active set $\overline{\mathcal{I}}$ or reduce the smoothing parameter $\alpha$ as shown in Algorithm~\ref{alg:Adjust-active-set}. The Newton steps then continue until either the active set is changed or the smoothing parameter is reduced. If the Newton method completes with $\alpha\leq\alpha_{\min}$, then the algorithm terminates. The parameter $0<\eta<1$ is used to identify completion of the Newton method for a given value of $\alpha$: $\left|\bm{d}^{T} \bm{g}\right|$ estimates the reduction of the objective function, which should be small compared to $\alpha$. We chose $\eta=1/10$ for our implementation.

\begin{algorithm}
\begin{algorithmic}[1]
    \State$\mathcal{J} \gets \left\{ \,j \in\mathcal{I}\mid\left|\widehat{g}_{j}\right| > \mu\,\right\} $
    \If{$|\mathcal{J}| > 0$ and $adjust$}
      \State$\mathcal{I} \gets \mathcal{I} \setminus \mathcal{J}$; \quad $\overline{\mathcal{I}} \gets \overline{\mathcal{I}} \cup \mathcal{J}'$; \Comment{Adjust active set}
      \State $g_j \gets g_j - \mu \sign(\widehat{g}_j) \quad \forall j \in \mathcal{J}$
      \State $\bm{d}_{\mathcal{J}} \gets -\bm{g}_{\mathcal{J}}$ \Comment{Gradient descent step}
      \State$adjust \gets \texttt{false}$ \Comment{Adjust active set once for each $\alpha$}
    \Else 
      \State $\alpha \gets \beta\,\alpha$ \Comment{Reduce $\alpha$ and optimize for this new $\alpha$}
      \State $adjust \gets \texttt{true}$
      \State $\mathbf{continue}$
    \EndIf
\end{algorithmic}
\caption{Adjust active set \& reduce smoothing parameter}
\label{alg:Adjust-active-set}
\end{algorithm}

Algorithm \ref{alg:Adjust-active-set} desribes when to adjust the active set or to reduce the smoothing parameter $\alpha$. If the smoothed problem is nearly solved for the current $\overline{\mathcal{I}}$ and $\alpha$, the set $\mathcal{J}$ of elements $j$ with $\left|\widehat{g}_{j}\right| > \mu$ is included to the active set and then a gradient descent step is taken. At this point, we take a gradient descent step rather than a Newton step in order to guarantee that $\bm{d}$ is a descent direction; that is, $\bm{d}^T \bm{g} < 0$. It is worth noting that the active set is adjusted once for each $\alpha$. The reason is to avoid adding and removing the same elements from the active set again and again. The Boolean variable $adjust$ indicates whether the active set would be adjusted in that iteration. 

Notice that elements can be added to $\mathcal{I}$ (line~20 of Algorithm~\ref{alg:smoothed-l1-Newton-step}) as well as removed from $\mathcal{I}$ (line~3 of Algorithm~\ref{alg:Adjust-active-set}). Note that multiple elements can be removed from $\overline{\mathcal{I}}$ in a single iteration. This helps keeping the size of active set small, which plays a significant role in the speed of the algorithm.

The line search algorithm is shown as Algorithm~\ref{alg:l1-linesearch}. Note that $0<c_{1}<1$ is the sufficient decrease parameter \cite[p.~33]{Noc2006}. The actual implementation differs slightly from the pseudo-code in that the recomputation of $\mathcal{I}$ on line~14 of Algorithm~\ref{alg:smoothed-l1-Newton-step} uses some additional information returned from \textsc{LinesearchL1}: In floating point arithmetic there is no guarantee that $\mathcal{I}\gets\left\{ \,j\mid w_{j}=0\,\right\} $ will identify components $w_{j}$ that would be set to zero in exact arithmetic. Specifically, setting $s\gets-w_{j}/d_{j}$ does not ensure that $w_{j}+s\,d_{j}$ evaluates to zero in floating point arithmetic. Therefore, the linesearch function \textsc{LinesearchL1} actually returns both $s$ and $j_{1}$ and $j_{2}$: if $j_{1}=j_{2}=j$, then $s=-w_{j}/d_{j}$ for $j=j_{1}=j_{2}$ and we would set $w_{j}+s\,d_{j}=0$ and the element $j$ is added to the inactive set $\mathcal{I}$.

\subsection{Issues with the line search} \label{subsec:Issues-line-search}
Line searches are often needed in optimization algorithms because the predicted step ``goes too far'', or fails to decrease the objective function value significantly. Typically, machine learning methods avoid line searches as each function evaluation requires a pass over a large dataset. This cost can be mitigated if we rarely need to use anything other than the initial step length. Let $\bm{d}$ denote the Newton direction. If the quadratic Taylor polynomial at $s=0$ for $f(\bm{w}+s\bm{d}; \mathcal{D})$ is a poor approximation to $f(\bm{w}+s\bm{d};\mathcal{D})$, then it may be necessary to perform many line search steps, which will require many function evaluations. As noted above in Section~\ref{subsec:Smoothed-hinge-loss}, for $N$ large, $N^{-1}\sum_{i=1}^{N}\varphi(\bm{x}_{i},y_{i};\bm{w})$ is close to a smooth function $h(\bm{w})$. The $\ell^{1}$\,SVM problem uses non-smooth $R(\bm{w})$, so $R(\bm{w})$ must be represented, at least approximately, in the line search procedure. Traditional Newton methods initially use the step length $s=1$, which will satisfy the line search criteria when $\bm{w}$ is close to the optimal $\bm{w}^{*}$, and the Armijo/backtracking line search (\cite{Arm66}, \cite[p.~33]{Noc2006}) requires few function evaluations while improving robustness of the algorithm. However, with non-smooth functions, such as the $\ell^{1}$ penalty, this choice can result in many function evaluations for a single line search. We therefore developed a line search method that incorporates the $\ell^1$ penalty.

\begin{algorithm}[H]
\begin{algorithmic}[1]
\Require{$a,\,\mu,\, s_{\max} \geq 0$, $\bm{d} \neq \bm{0}$, and either $a>0$ or $\mu\left\Vert \bm{d}\right\Vert _{1}>-b$}

\Function{LinesearchL1}{$\bm{w},\bm{d},b,a,\mu,s_{\max}$}

\Comment{It returns $s$ that minimizes $as^{2}+bs+\mu\left\Vert \bm{w}+s\bm{d}\right\Vert _{1}$
over $0\leq s\leq s_{\max}$}
  \State$m \gets \text{dimension}(\bm{w})$
  \State{}find function $p\colon\left\{ 1,2,\ldots,m\right\} \to\left\{ 1,2,\ldots,m\right\}$ such that
  \State$\qquad$$\text{range}(p)=\left\{ \,j\mid-w_{j}/d_{j}>0\,\right\}$ and $[-w_{p(i)}/d_{p(i)}]_{i=1}^{m}$ is sorted
  \State$i_{1}\gets0$; $j_{1}\gets0$; $s_{1}\gets0$; $\mathit{slope}_{1}=b+\mu\,\text{sign}(\bm{w})^{T}\bm{d}$
  \State$i_{2}\gets m+1$; $j_{2}\gets m+1$; $s_{2}\gets+\infty$;
$\mathit{slope}_{2}\gets\begin{cases}
+\infty, & \text{if }a>0,\\
b+\mu\left\Vert \bm{d}\right\Vert _{1}, & \text{if }a=0.
\end{cases}$
  \State\textbf{if }$\mathit{slope}_{1}\geq0$ \textbf{then} $\mathbf{return}\;s_{1}$;
  \State\textbf{if }$\mathit{slope}_{2}\leq0$ \textbf{then} $\mathbf{return}\;s_{2}$;
\While{$i_{2}>i_{1}+1$}\Comment{binary search}
  \State$i\gets\left\lfloor (i_{1}+i_{2})/2\right\rfloor $; $j\gets p(i)$
  \State$s\gets-w_{j}/d_{j}$\Comment{compute slopes on either side of $s$}
  \State$\mathit{slope}_{0}\gets2as+b+\sum_{k:k\neq j}\text{sign}(w_{k}+sd_{k})d_{k}$
  \State$\mathit{slope}_{+}\gets\mathit{slope}_{0}+\mu\left|d_{j}\right|$;
$\mathit{slope}_{-}\gets\mathit{slope}_{0}-\mu\left|d_{j}\right|$
  \If{($\mathit{slope}_{-}=0$ or $\mathit{slope}_{+}=0$) or ($\mathit{slope}_{-}<0$ and $\mathit{slope}_{+}>0$)}
    \State$\mathbf{return}\;s$
  \Else{}  \If{ $\mathit{slope}_{+}<0$}
      \State$i_{1}\gets i$; $j_{1}\gets j$; $\mathit{slope}_{1}\gets\mathit{slope}_{-}$; $s_{1}\gets s$
    \Else{}
    \State$i_{2}\gets i$; $j_{2}\gets j$; $\mathit{slope}_{2}\gets\mathit{slope}_{+}$; $s_{2}\gets s$  \EndIf
  \EndIf
\EndWhile
\Comment{Note that $i_{2}=i_{1}+1$ \& the optimal $s$ is in $(s_{1},s_{2})$}
\State$s\gets(s_{1}\mathit{slope}_{2}-s_{2}\mathit{slope}_{1})/(\mathit{slope}_{2}-\mathit{slope}_{1})$

\EndFunction
\end{algorithmic}
\caption{\label{alg:l1-linesearch}Line search algorithm for quadratic plus $\ell^{1}$ penalty}
\end{algorithm}

Including the $\ell^{1}$ penalty, we consider the minimization problem
\[
\min_{s \geq 0}\ \widehat{f}_\alpha(\bm{w}+s\bm{d})+\mu\left\Vert \bm{w}+s\bm{d}\right\Vert_{1},
\]
where $\widehat{f}_\alpha$ is a smooth function. Since we can estimate the Hessian matrices accurately, we can use a quadratic approximation for $\widehat{f}_\alpha(\bm{w}+s\bm{d})\approx a\,s^{2}+b\,s+c$, where
\[
a := \frac{1}{2} \bm{d}^T H \bm{d}, \qquad b := \bm{d}^T \widehat{\bm{g}}, \qquad c := \widehat{f}_{\alpha}(\bm{w}).
\]
Our line search then seeks to minimize 
\begin{equation} \label{eq:line-search-quad}
j(s) := a\,s^{2}+b\,s+c+\mu\left\Vert \bm{w}+s\bm{d}\right\Vert_{1}  \qquad \text{over } s \geq 0.
\end{equation}
We assume that $\bm{d}$ is a descent direction, so that $j'(0)<0$. Provided $a,\,\mu\geq0$, this is a convex function, and hence the derivative $j'(s)$ is a non-decreasing function of $s$. Provided $\left\Vert \bm{d}\right\Vert_{1}>b$ \emph{or} $a>0$ \emph{or} $\mu>0$, there is a global minimizer of $j$; if $a>0$ then it is unique. The task is to compute this minimizer efficiently. The minimizer is characterized by either $j'(s)=0$, or 
\[
j'(s^{-}):=\lim_{r\to s^{-}}j'(r)\leq0 \quad \text{and} \quad j'(s^{+}):=\lim_{r\to s^{+}}j'(r)\geq0.
\]
This can be done using a binary search algorithm that just uses $a$, $b$, $\mu$, $\bm{w}$ and $\bm{d}$, \emph{and no additional evaluations of} $\theta$. Note that 
\begin{equation} \label{eq:j-prime}
j'(s) = 2as + b + \mu \sum_{i=1}^{m} \text{sign} (w_{i}+sd_{i})\,d_{i}.
\end{equation}
If $\mu=0$ and $a>0$, then the minimizer is clearly $s=-b/(2a)$. Assuming $a>0$ and $\mu\geq0$, the minimizing value of $s$ must lie in the interval $(0,\,s_{\max}]$ where $s_{\max}=(\left|b\right|+\mu\left\Vert \bm{d}\right\Vert_{1})/(2a)$. The points of discontinuity of $j'(s)$ are $\sigma_{i} = -w_{i}/d_{i}$, $i=1,\,2,\,\ldots,\,m$. If $d_{i}=0$, we simply ignore $\sigma_{i}$. Let $\left\{ \widehat{\sigma}_{1},\,\widehat{\sigma}_{2},\,\ldots,\,\widehat{\sigma}_{r}\right\} =\left\{ \sigma_{i}\mid\sigma_{i}>0\right\} $ with $\widehat{\sigma}_{1}<\widehat{\sigma}_{2}<\cdots<\widehat{\sigma}_{r}$. Set $\widehat{\sigma}_{0}=0$. We check $j'(\widehat{\sigma}_{0}^{+})$ and $j'(\widehat{\sigma}_{r}^{+})$:
\begin{itemize}
\item If $j'(\widehat{\sigma}_{0}^{+})\geq0$, then the optimal $s$ is $s^{*}=0=\widehat{\sigma}_{0}$.
\item If $a=0$ and $j'(\widehat{\sigma}_{r}^{+})<0$, then $j(s)\to-\infty$ as $s\to\infty$ and there is no minimum.
\item If $a>0$ and $j'(\widehat{\sigma}_{r}^{+})<0$, then the optimal $s$ is
\[
s^{*} = -\frac{1}{2a}\left(b+\mu\sum_{i=1}^{m}\text{sign}(d_{i})d_{i}\right)
= -\frac{b+\mu\left\Vert \bm{d}\right\Vert_{1}}{2a}
= \widehat{\sigma}_{r} - \frac{j'(\widehat{\sigma}_{r}^{+})}{2a}
> \widehat{\sigma}_{r},
\]
since $j'(\widehat{\sigma}_{r}^{+}) = 2a\widehat{\sigma}_{r} + b + \mu\left\Vert \bm{d}\right\Vert_{1}$.
\end{itemize}
Now consider the sequence 
\begin{equation} \label{eq:j-prime-sequence}
j'(\widehat{\sigma}_{0}^{+}),\,j'(\widehat{\sigma}_{1}^{-}),\,j'(\widehat{\sigma}_{1}^{+}),\,j'(\widehat{\sigma}_{2}^{-}),\,j'(\widehat{\sigma}_{2}^{+}),\,\ldots,\,j'(\widehat{\sigma}_{r}^{+}).
\end{equation}
Since this sequence is a non-decreasing sequence, if $j'(\widehat{\sigma}_{0}^{+})<0$ and $j'(\widehat{\sigma}_{r}^{+})>0$, then it crosses from being $\leq0$ to $>0$ at some point. If $j'(\widehat{\sigma}_{i}^{\pm})=0$ for some $i$ and choice of sign, then $s^{*}=\widehat{\sigma}_{i}$. Thus, we assume without loss of generality that $j'(\widehat{\sigma}_{i}^{\pm})\neq0$ for any $i$ and choice of sign. In this case, there exists an index $i$ such that \emph{either}
\begin{itemize}
\item[(i)] $j'(\widehat{\sigma}_{i}^{-})<0$ and $j'(\widehat{\sigma}_{i}^{+})>0$, \emph{or}
\item[(ii)] $j'(\widehat{\sigma}_{i}^{+})<0$ and $j'(\widehat{\sigma}_{i+1}^{-})>0$.
\end{itemize}
In the former case, we get $s^{*}= \widehat{\sigma}_{i}$. In the latter case, for $s\in(\widehat{\sigma}_{i},\widehat{\sigma}_{i+1})$, we have $j'(s)=j'(\widehat{\sigma}_{i}^{+})+2a(s-\widehat{\sigma}_{i})$, and hence $s^{*}=\widehat{\sigma}_{i}-j'(\widehat{\sigma}_{i}^{+})/(2a)$.

Finding the point where the sequence (\ref{eq:j-prime-sequence}) crosses zero can be carried out by binary search. Thus, it can be computed in $\mathcal{O}(\log m)$ time as $r\leq m$ after sorting the positive $\sigma_i$'s which takes $\mathcal{O}(m \log m)$. A complete line search algorithm is outlined in Algorithm~\ref{alg:l1-linesearch}.

If the optimal value for $s^{*}$ is zero, then $\bm{d}$ is not a descent direction and so some other direction should be used. This can only occur if $\sigma_{i}=0$ for some $i$, indicating that $w_{i}=0$. Then in this case, we need to remove $w_{i}$ from the set of active variables.

\section{Convergence rate} \label{sec:convergence-rate}
The unresolved question about this algorithm is its efficiency. Standard results on Newton's method with line search applied to strictly convex functions indicate that the method will converge. However, the question we wish to answer is: \emph{How quickly will the method converge?} Ideally we would like a constant number of iterations of the Newton method between steps where $\alpha$ is reduced by a constant factor $\alpha\gets\beta\alpha$. This would lead to $\mathcal{O}(\log(1/\alpha_{\min}))$ gradient and Hessian evaluations to achieve an error tolerance of $\alpha_{\min}>0$. Furthermore, this number should be independent of $N$ and $\alpha>0$, and perhaps even of the dataset $\left\{ \,(\bm{x}_{i},y_{i})\mid i=1,\ldots,N\,\right\} $. In this section we do not discuss the active set aspects of the algorithm. We simply note that the $\ell^{1}$ penalty keeps the dimension small to moderate, and that active set changes only occur near the minimum for a given active set (when $\left|\bm{d}^{T}{\bm{g}}\right|<\eta\alpha$). Both of these aspects should keep the active sets from changing often, and slowing convergence \cite[pp.~265--275]{Fle87}.

We do not claim to be able to prove the rate of convergence, with asymptotically $\mathcal{O}(1)$ Newton steps for each reduction in the smoothing parameter $\alpha$. However, we give a partial argument for why we expect good performance for this algorithm. We separate the steps into three phases:
\begin{enumerate}
\item \emph{Opening}: the Newton steps with $\alpha\gg1$
\item \emph{Midgame}: the steps with $1\gg\alpha\gg1/N$
\item \emph{Endgame}: the steps with $1/N\gg\alpha>0$
\end{enumerate}
The \emph{opening} can be treated in a straightforward asymptotic approach; not much accuracy is expected, but the computed weight vector is ``in the right ball-park''. The \emph{midgame} is more difficult to analyze, but with $N\alpha\gg1$, assuming the data points $(\bm{x}_{i},y_{i})$ are picked independently from a common probability distribution with a Lipschitz probability density function, we can use expectations and variances to estimate Hessian matrices and other quantities to show that Newton's method with line search converges rapidly and reliably. The \emph{endgame} analysis assumes that we have separated out the relevant data points for the support vector and assumes that some non-degeneracy properties hold. This analysis is based on self-similarity of the smoothed optimization problems as $\alpha$ decreases.

Between these asymptotic regimes are asymptotic regimes $\alpha=\Theta(1)$ and $\alpha=\Theta(N^{-1})$. While these asymptotic bands are transitions between regions covered by our analysis, and so far escape our understanding of how the algorithm behaves in these regions. Nevertheless, we believe that the algorithm can successfully traverse these transitions with only a modest number of guarded Newton steps. More discussion of these
issues can be found in Section~\ref{subsec:Summary-of-convergence}.

In our analysis we assume that the data points $(\bm{x}_{i},y_{i})$ are drawn independently from a fixed probability distribution with a Lipschitz density function $p(\bm{x},y)$. We use $\mathcal{D}$ to denote the dataset drawn from this distribution and $\E[Z]$ to denote expectation of a random variable $Z$ with respect to the generation of the dataset in this fashion, while $\Var[\bm{Z}]$ is used to denote the variance\textendash covariance matrix of $\bm{Z}$ for vector-valued $\bm{Z}$. For matrix-valued random variables $Z$, we use $\Var\left[Z\right]$ to denote the scalar quantity $\E\left[\left(Z-\E\left[Z\right]\right)\bullet\left(Z-\E\left[Z\right]\right)\right]$ where $A\bullet B=\text{trace}(A^{T}B)$ is the Frobenius inner product of matrices $A$ and $B$ \cite[p.~332, Ex.~24]{Hor85}. Note that the usual properties of variances hold for this variance: $\Var\left[a\,Z\right]=a^{2}\Var\left[Z\right]$; $\Var\left[Z_{1}+Z_{2}\right]=\Var\left[Z_{1}\right]+\Var\left[Z_{2}\right]$ for independent $Z_{1}$ and $Z_{2}$; and 
\[
\Var\left[Z\right]=\E\left[Z\bullet Z\right]-\E\left[Z\right]\bullet\E\left[Z\right]=\E\left[\left\Vert Z\right\Vert _{F}^{2}\right]-\left\Vert \E\left[Z\right]\right\Vert _{F}^{2}.
\]
The expectation $\E[\widehat{f}_{\alpha}(\bm{w};\mathcal{D})]$ is the expectation of $\widehat{f}_{\alpha}(\bm{w};\mathcal{D})$ over all datasets $\mathcal{D}$ with each $(\bm{x}_{i},y_{i})$, $i=1,2,\ldots,N$ drawn independently from the probability distribution with density $p(\bm{x},y)$. That is,
\begin{align*}
\E[\widehat{f}_{\alpha}(\bm{w};\mathcal{D})] & = \E_{(\bm{X},Y)\sim p}\left[\alpha\,\psi\left(\frac{1-Y\bm{X}^{T}\bm{w}}{\alpha}\right)\right]+\frac{1}{2}\lambda\left\Vert \bm{w}\right\Vert _{2}^{2}\\
&=\sum_{y\in\left\{ \pm1\right\} }\int_{\R^{m}}p(\bm{x},y)\,\alpha\,\psi\left(\frac{1-y\bm{x}^{T}\bm{w}}{\alpha}\right)\,d\bm{x}+\frac{1}{2}\lambda\left\Vert \bm{w}\right\Vert_{2}^{2}.
\end{align*}
By linearity of differentiation and expectations,
\[
\nabla\,\E[\widehat{f}_{\alpha}(\bm{w};\mathcal{D})]=\E[\nabla\,\widehat{f}_{\alpha}(\bm{w};\mathcal{D})] \ \text{ and }\ \Hess \E[\widehat{f}_{\alpha}(\bm{w};\mathcal{D})]=\E[\Hess \widehat{f}_{\alpha}(\bm{w};\mathcal{D})].
\]
Note that $\widehat{f}_{\alpha}(\bm{w}; \mathcal{D})$ is the expected value of $\alpha\,\psi((1-Y\bm{X}^{T}\bm{w})/\alpha)+\frac{1}{2}\lambda\left\Vert \bm{w}\right\Vert _{2}^{2}$ where $(\bm{X},Y)$ are distributed according to the points $(\bm{x}_{i},y_{i})$ with equal probability $1/N$. Certain results are independent of the probability distribution of the data points (for example, $\bm{w}\mapsto\E[\Hess \widehat{f}_{\alpha}(\bm{w}; \mathcal{D})]$ is Lipschitz with constant $\mathcal{O}(1/\alpha^{2})$), while certain results require the assumption that the probability density $\bm{x}\mapsto p(\bm{x},y)$ is Lipschitz. The integration techniques used in the results obtained below involve the marginal probability distribution:
\[
p_{\bm{w}}(s,y)=\int_{\left\{ \bm{w}\right\} ^{\perp}}p(s\bm{w}+\bm{v})\,d\bm{v}.
\]
It should also be noted that the asymptotic results are said to hold ``with probability one''. The justification for these claims is that for any family of non-negative random variables $Z_{\alpha}$, if $\E\left[Z_{\alpha}\right]\leq h(\alpha)$ then by Chebyshev's inequality, $\Pr\left[Z_{\alpha}>C\,h(\alpha)\right]\leq1-1/C$ for any positive $C$. Thus $Z_{\alpha}=\mathcal{O}(h(\alpha))$ with probability one.

\subsection{Assumptions} \label{subsec:Assumptions}
It is assumed that the data points $(\bm{x}_{i},y_{i})$, $i=1,2,\ldots,N$ are taken independently from a distribution with probability density function $p(\bm{x},y)$ ($\bm{x}\in\R^{m}$ and $y\in\left\{ \pm1\right\} $) that is Lipschitz continuous in $\bm{x}$ and has bounded support. In what follows, $(\bm{X},Y)$ is a random vector with this probability distribution.

Another issue is that the solution $\bm{w}=\bm{0}$ can minimize $\widehat{f}(\bm{w}; \mathcal{D})$, but should be rare. In this case we have $\nabla\widehat{f}(\bm{w}; \mathcal{D})=\bm{0}$, and hence
\[
\sum_{i=1}^{N}y_{i}\bm{x}_{i}=\bm{0}.
\]
The case where $\bm{w}=\bm{0}$ means that the SVM cannot distinguish between any data points. In order for the SVM to distinguish any data point $(\bm{x}_{i},y_{i})$, we need to have $1<y_{i}\bm{x}_{i}^{T}\bm{w}$ so that $\left\Vert \bm{w}\right\Vert_{2}>1/\left\Vert \bm{x}_{i}\right\Vert_{2}\geq1/\max_{i}\left\Vert \bm{x}_{i}\right\Vert_{2}$.

We assume that $p(\bm{x},y)=0$ whenever $\left\Vert \bm{x}\right\Vert_{2}\geq R$ for a given $R \in \R$. Thus, for any dataset $\mathcal{D}$ drawn from this distribution, we have $\max_{i}\left\Vert \bm{x}_{i}\right\Vert_{2}\leq R$, and so any useful optimal $\bm{w}$ has $\left\Vert \bm{w}\right\Vert_{2}>1/R$. 

We assume that either $\lambda>0$ or that the optimal $\bm{w}=\bm{w}^{*}$ for minimizing $\E\big[\widehat{f}(\bm{w};\mathcal{D})\big]$ has the property that 
\[
\sum_{y\in\left\{ \pm1\right\} }\int_{\left\{ \bm{x}\mid y=\bm{x}^{T}\bm{w}\right\} }p(\bm{x},y)\,dS(\bm{x})>0.
\]
Note that either of these conditions ensures that $\Hess \E\big[\widehat{f}(\bm{w}; \mathcal{D})\big]$ is positive definite at the optimal $\bm{w}=\bm{w}^{*}$. 

\subsection{Opening} \label{subsec:Opening}
We begin with $\alpha>0$ large. Note that as $\alpha\to\infty$, we have
\begin{align*}
\Hess \widehat{f}_{\alpha}(\bm{w}; \mathcal{D}) &= \frac{1}{N}\sum_{i=1}^{N}\frac{1}{\alpha}\psi'' \left(\frac{1-y_{i}\bm{w}^{T}\bm{x}_{i}}{\alpha} \right) \bm{x}_{i}\bm{x}_{i}^{T} + \lambda I\to\lambda I, \\
\nabla\widehat{f}_{\alpha}(\bm{w}; \mathcal{D}) &= \frac{1}{N}\sum_{i=1}^{N}\psi' \left(\frac{1-y_{i}\bm{w}^{T}\bm{x}_{i}}{\alpha}\right)(-y_{i}\bm{x}_{i}) + \lambda\bm{w} \to -\frac{\psi'(0)}{N} \sum_{i=1}^{N}y_{i}\bm{x}_{i} +\lambda\bm{w},
\end{align*}
and we are in the asymptotic regime where $\alpha\gg1+\left\Vert \bm{w}\right\Vert \max_{i}\left\Vert \bm{x}_{i}\right\Vert$. Thus, we have
\[
\bm{w}-\big(\Hess \widehat{f}_{\alpha}(\bm{w}; \mathcal{D}) \big)^{-1}\nabla\widehat{f}_{\alpha}(\bm{w}; \mathcal{D})\to\frac{1}{2N\lambda}\sum_{i=1}^{N}y_{i}\bm{x}_{i}\qquad\text{as }\alpha\to\infty,
\]
since $\psi'(0)=\frac{1}{2}$. Thus, for large $\alpha$, one step of Newton's method will bring the value of $\bm{w}$ within $\mathcal{O}(1)$ of the optimal, provided $\lambda>0$.

\subsection{Midgame} \label{subsec:Midgame}
For the midgame estimates, we assume that $1\gg\alpha\gg1/N$ so that we can use averages with respect to the underlying probability distribution
\[
{\displaystyle \sum_{y\in\left\{ \pm1\right\} }\int_{\R^{m}}}p(\bm{x},y)\,\varphi_{\alpha}(\bm{x},y;\bm{w})\,d\bm{x}
\]
in place of sums over the dataset: $(1/N)\sum_{i=1}^{N}\varphi_{\alpha}(\bm{x}_{i},y_{i};\bm{w})$. The aim is to bound the number of Newton iterations each time the smoothing parameter is reduced: $\alpha\gets\beta\alpha$. The results of Boyd, Boyd, and Vandenberghe \cite[Chap.~9]{Boy2004} are used to do this; the bounds are based on several quantities: upper and lower bounds on the eigenvalues of $\Hess \widehat{f}_{\alpha}(\bm{w}; \mathcal{D})$, Lipschitz constants for $\bm{w}\mapsto\Hess \widehat{f}_{\alpha}(\bm{w}; \mathcal{D})$, and the parameters of the Armijo backtracking line search used.

The basis for the results of this section is the Lipschitz continuity of the probability density function $p(\bm{x},y)$. Since the data items of $\mathcal{D}=\left\{ \,(\bm{x}_{i},y_{i})\mid i=1,2,\ldots,N\,\right\} $ are sampled independently from this probability distribution, we need to consider both the expected values of crucial quantities and the variances of these quantities. Section~\ref{subsec:Lipschitz-EHessf} gives a Lipschitz constant for $\E[\Hess \widehat{f}_{\alpha}(\bm{w};\mathcal{D})]$. Section~\ref{subsec:Lips-const-Hessf} gives a asymptotic estimates for the Lipschitz constants for $\Hess \widehat{f}_{\alpha}(\bm{w};\mathcal{D})$ using estimates for $\Var\big[\Hess \widehat{f}_{\alpha}(\bm{w};\mathcal{D})\big]$. Section~\ref{subsec:Uniform-bdd-condEHess} gives upper and lower bounds on the eigenvalues of $\E\big[\Hess \widehat{f}_{\alpha}(\bm{w};\mathcal{D})\big]$ in the form of bounds $aI\preceq\E\big[\Hess \widehat{f}_{\alpha}(\bm{w};\mathcal{D})\big]\preceq bI$. Note that $A\preceq B$ means that $\bm{z}^{T}A\bm{z} \leq \bm{z}^{T}B\bm{z}$ for all $\bm{z}$. Section~\ref{subsec:Var-Hess} gives bounds on $\Var\big[\Hess \widehat{f}_{\alpha}(\bm{w};\mathcal{D})\big]$, so that we have asymptotic estimates for eigenvalue of $\Hess \widehat{f}_{\alpha}(\bm{w};\mathcal{D})$ with probability one. Section~\ref{subsec:Variance-of-gradients} establishes bounds for $\Var\big[\nabla\widehat{f}_{\alpha}(\bm{w};\mathcal{D})\big]$. This is combined with the next section, Section~\ref{subsec:Changing-alpha} which gives bounds on $\big\Vert \E\big[\nabla\widehat{f}_{\beta\alpha}(\bm{w};\mathcal{D})\big]\big\Vert $ so that we can estimate the gradient and function values after reducing alpha: $\alpha\gets\beta\alpha$. Section~\ref{subsec:midgame-convce-Newton} combines the results of the previous sections to give the result that with probability one, $\mathcal{O}(1+\log\log(1/\alpha))$ Newton iterations are needed with each step if $1\gg\alpha\gg N^{-1}$. The section also shows that only $\mathcal{O}(1+\log\log(1/\alpha))$ function values are needed in the line search. 

\subsubsection{Lipschitz continuity of $\protect\E[\Hess \widehat{f}_{\alpha}(\bm{w};\mathcal{D})]$ in $\bm{w}$} \label{subsec:Lipschitz-EHessf}
We wish to show a Lipschitz continuity property of $\E[\Hess \widehat{f}_{\alpha}(\bm{w};\mathcal{D})]$ for $\left\Vert \bm{w}\right\Vert \geq1/(2R)$. 
\begin{lemma}
Assuming that $\bm{x}\mapsto p(\bm{x},y)$ is Lipschitz, for any compact set $C$ not containing $\bm{0}$, there is a constant $L$ independent of $\alpha$ such that
\[
\left\Vert \E[\Hess \widehat{f}_{\alpha}(\bm{w};\mathcal{D})]-\E[\Hess \widehat{f}_{\alpha}(\bm{w}';\mathcal{D})]\right\Vert_{2}\leq L\left\Vert \bm{w}-\bm{w}'\right\Vert, \quad \forall\, \bm{w},\,\bm{w}'\in C.
\]
\end{lemma}
\begin{proof}
Let $A=I+(\bm{w}'-\bm{w})\bm{w}^{T}/\left\Vert \bm{w}\right\Vert _{2}^{2}$
so $A\bm{w}=\bm{w}'$. We can see that $(\bm{w}')^{T}\bm{x}=(A\bm{w})^{T}\bm{x}=\bm{w}^{T}(A^{T}\bm{x})$.
Putting $\bm{x}'=A^{T}\bm{x}$, we get
\[
\left\Vert \bm{x}'-\bm{x}\right\Vert _{2}=\frac{\left\Vert \bm{w}(\bm{w}'-\bm{w})^{T}\bm{x}\right\Vert _{2}}{\left\Vert \bm{w}\right\Vert _{2}^{2}}\leq\frac{\left\Vert \bm{w}'-\bm{w}\right\Vert _{2}}{\left\Vert \bm{w}\right\Vert _{2}}\left\Vert \bm{x}\right\Vert _{2}.
\]
Using the expectation approximation (valid for $N\alpha\gg1$), 
\begin{align*}
 & \E\left[\Hess\widehat{f}_{\alpha}(\bm{w}';\mathcal{D})\right]-\lambda I\\
 & \qquad{}=\sum_{y\in\left\{ \pm1\right\} }\int_{\R^{m}}p(\bm{x},y)\,\frac{1}{\alpha}\psi''\left(\frac{1-y(\bm{w}')^{T}\bm{x}}{\alpha}\right)\bm{x}\bm{x}^{T}d\bm{x}\\
 & \qquad{}=\sum_{y\in\left\{ \pm1\right\} }\int_{\R^{m}}p(\bm{x},y)\,\frac{1}{\alpha}\psi''\left(\frac{1-y\bm{w}^{T}A^{T}\bm{x}}{\alpha}\right)\bm{x}\bm{x}^{T}d\bm{x}\\
 & \qquad{}=\sum_{y\in\left\{ \pm1\right\} }\int_{\R^{m}}p(A^{-T}\bm{x}',y)\,\frac{1}{\alpha}\psi''\left(\frac{1-y\bm{w}^{T}\bm{x}'}{\alpha}\right)\left|\det(A^{-T})\right|\,(A^{-T}\bm{x}')(A^{-T}\bm{x}')^{T}d\bm{x}'\\
 & \qquad{}=\left|\det(A^{-T})\right|A^{-T}\sum_{y\in\left\{ \pm1\right\} }\int_{\R^{m}}p(A^{-T}\bm{x}',y)\,\frac{1}{\alpha}\psi''\left(\frac{1-y\bm{w}^{T}\bm{x}'}{\alpha}\right)\,\bm{x}'(\bm{x}')^{T}d\bm{x}'\,A^{-1};\\
 & \E\left[\Hess\widehat{f}_{\alpha}(\bm{w}';\mathcal{D})-\Hess\widehat{f}_{\alpha}(\bm{w};\mathcal{D})\right]\\
 & \qquad{}=\E\left[\Hess\widehat{f}_{\alpha}(\bm{w}';\mathcal{D})-\left|\det(A^{-T})\right|A^{-T}\Hess\widehat{f}_{\alpha}(\bm{w};\mathcal{D})A^{-1}\right]\\
 & \qquad\qquad{}+\E\left[\left|\det(A^{-T})\right|A^{-T}\Hess\widehat{f}_{\alpha}(\bm{w};\mathcal{D})A^{-1}-\Hess\widehat{f}_{\alpha}(\bm{w};\mathcal{D})\right].
\end{align*}
For the first part, we have
\begin{align*}
 & \left\Vert \E\left[\Hess\widehat{f}_{\alpha}(\bm{w}';\mathcal{D})-\left|\det(A^{-T})\right|A^{-T}\Hess\widehat{f}_{\alpha}(\bm{w};\mathcal{D})A^{-1}\right]\right\Vert _{2}\\
 & \qquad{}=\left|\det(A^{-T})\right|\left\Vert A^{-T}\sum_{y\in\left\{ \pm1\right\} }\int_{\R^{m}}\left[p(A^{-T}\bm{x},y)-p(\bm{x},y)\right]\times{}\right.\,\\
 & \qquad\qquad\qquad\qquad\qquad\left.{}\times\frac{1}{\alpha}\psi''\left(\frac{1-y\bm{w}^{T}\bm{x}}{\alpha}\right)\,\bm{x}\bm{x}^{T}\,d\bm{x}\,A^{-1}\right\Vert _{2}\\
 & \qquad{}\leq\left|\det(A^{-T})\right|\left\Vert A^{-1}\right\Vert _{2}^{2}\sum_{y\in\left\{ \pm1\right\} }\int_{\left\{ \bm{x}\mid\left\Vert \bm{x}\right\Vert _{2}\text{ or }\left\Vert A^{-1}\bm{x}\right\Vert _{2}\leq R\right\} }\left|p(A^{-T}\bm{x},y)-p(\bm{x},y)\right|\times{}\\
 & \qquad\qquad\qquad\qquad\qquad{}\times\frac{1}{\alpha}\psi''\left(\frac{1-y\bm{w}^{T}\bm{x}}{\alpha}\right)\,\left\Vert \bm{x}\bm{x}^{T}\right\Vert _{2}\,d\bm{x}\qquad(\text{as }\left\Vert A^{-1}\right\Vert _{2}=\left\Vert A^{-T}\right\Vert _{2})\\
 & \qquad{}\leq\left|\det(A^{-T})\right|\left\Vert A^{-1}\right\Vert _{2}^{2}\sum_{y\in\left\{ \pm1\right\} }\\
 & \qquad\qquad\int_{\left\{ \bm{x}\mid\left\Vert \bm{x}\right\Vert _{2}\text{ or }\left\Vert A^{-T}\bm{x}\right\Vert _{2}\leq R\right\} }L_{p}\left\Vert A^{-T}\bm{x}-\bm{x}\right\Vert _{2}\,\frac{1}{\alpha}\psi''\left(\frac{1-y\bm{w}^{T}\bm{x}}{\alpha}\right)\,\left\Vert \bm{x}\right\Vert _{2}^{2}\,d\bm{x}\\
 & \qquad\qquad\qquad(\text{since }p\text{ is Lipschitz with constant }L_{p})\\
 & \qquad{}\leq\left|\det(A^{-T})\right|\left\Vert A^{-1}\right\Vert _{2}^{2}L_{p}\left\Vert A^{-T}-I\right\Vert _{2}\sum_{y\in\left\{ \pm1\right\} }\\
 & \qquad\qquad\int_{B_{R\max(1,\left\Vert A\right\Vert _{2})}}\frac{1}{\alpha}\psi''\left(\frac{1-y\bm{w}^{T}\bm{x}}{\alpha}\right)\,\left\Vert \bm{x}\right\Vert _{2}^{3}\,d\bm{x}\qquad(\text{since }\text{supp}\,p(\cdot,y)\subseteq B_{R})
\end{align*}

\begin{align*}
 & \qquad{}\leq\left|\det(A^{-T})\right|\left\Vert A^{-1}\right\Vert _{2}^{2}L_{p}\left\Vert A^{-T}-I\right\Vert _{2}\sum_{y\in\left\{ \pm1\right\} }\\
 & \qquad\qquad\int_{-\infty}^{+\infty}\int_{\left\{ \bm{w}\right\} ^{\perp}\cap B_{R\max(1,\left\Vert A\right\Vert _{2})}}\frac{1}{\alpha}\psi''\left(\frac{1-y\bm{w}^{T}(s\widehat{\bm{w}}+\bm{v})}{\alpha}\right)\,\left\Vert s\widehat{\bm{w}}+\bm{v}\right\Vert _{2}^{3}\,d\bm{v}\,ds\\
 & \qquad\qquad\qquad(\bm{x}=s\widehat{\bm{w}}+\bm{v}\text{ where }\widehat{\bm{w}}=\bm{w}/\left\Vert \bm{w}\right\Vert _{2}\text{ and }\bm{v}\perp\bm{w})\\
 & \qquad{}\leq\left|\det(A^{-T})\right|\left\Vert A^{-1}\right\Vert _{2}^{2}L_{p}\left\Vert A^{-T}-I\right\Vert _{2}\text{vol}_{m-1}\left(\left\{ \bm{w}\right\} ^{\perp}\cap B_{R\max(1,\left\Vert A\right\Vert _{2})}\right)\\
 & \qquad\qquad\qquad\,\left(R\max(1,\left\Vert A\right\Vert _{2})\right)^{3}\int_{-\infty}^{+\infty}\frac{1}{\alpha}\psi''\left(\frac{1-ys\left\Vert \bm{w}\right\Vert _{2}}{\alpha}\right)ds\\
 & \qquad{}\leq\left|\det(A^{-T})\right|\left\Vert A^{-1}\right\Vert _{2}^{2}L_{p}\left\Vert A^{-T}-I\right\Vert _{2}\text{vol}_{m-1}\left(\left\{ \bm{w}\right\} ^{\perp}\cap B_{R\max(1,\left\Vert A\right\Vert _{2})}\right)\\
 & \qquad\qquad\qquad\,\left(R\max(1,\left\Vert A\right\Vert _{2})\right)^{3}/\left\Vert \bm{w}\right\Vert _{2}.
\end{align*}
Note that this bound is not sharp; better constants can be found by more careful analysis of the behavior of $p(\bm{x},y)$. We have independence of $\alpha$ in the right-hand side since 
\[
\int_{-\infty}^{+\infty}\frac{1}{\alpha}\psi''\left(\frac{1-\omega s}{\alpha}\right)\,ds=\int_{-\infty}^{+\infty}\frac{\psi''(u)}{\left|\omega\right|}\,du=\frac{\lim_{u\to+\infty}\psi'(u)-\lim_{u\to-\infty}\psi'(u)}{\left|\omega\right|}=\frac{1}{\left|\omega\right|}.
\]
Since $(I+\bm{a}\bm{b}^{T})^{-1}=I-\bm{a}\bm{b}^{T}/(1+\bm{b}^{T}\bm{a})$
we have
\begin{align*}
\left\Vert A^{-T}-I\right\Vert _{2} & =\left\Vert \bm{w}(\bm{w}'-\bm{w})^{T}\right\Vert _{2}/(\left\Vert \bm{w}\right\Vert ^{2}\left[1+(\bm{w}'-\bm{w})^{T}\bm{w}\right])\\
 & \leq\left\Vert \bm{w}'-\bm{w}\right\Vert /\left(\left\Vert \bm{w}\right\Vert \left[1+(\bm{w}'-\bm{w})^{T}\bm{w}\right]\right).
\end{align*}
Thus, for $\bm{w}$ bounded away from zero and $\left\Vert \bm{w}'-\bm{w}\right\Vert \leq1/(2\left\Vert \bm{w}\right\Vert )$,
we see that 
\begin{align*}
\left\Vert A^{-T}-I\right\Vert _{2} \leq\left\Vert \bm{w}'-\bm{w}\right\Vert /\left(\left\Vert \bm{w}\right\Vert \left[1+(\bm{w}'-\bm{w})^{T}\bm{w}\right]\right)
\leq2\left\Vert \bm{w}'-\bm{w}\right\Vert /\left\Vert \bm{w}\right\Vert .
\end{align*}
For the second part,
\begin{align*}
\E\left[\left|\det(A^{-T})\right|A^{-T}\Hess\widehat{f}_{\alpha}(\bm{w};\mathcal{D})A^{-1}-\Hess\widehat{f}_{\alpha}(\bm{w};\mathcal{D})\right]
 =\mathcal{O}\left(\left\Vert \bm{w}'-\bm{w}\right\Vert /\left\Vert \bm{w}\right\Vert \right)
\end{align*}
since 
\begin{align*}
\det(A)-1 & =\det(I+(\bm{w}'-\bm{w})\bm{w}^{T}/\left\Vert \bm{w}\right\Vert ^{2})-1\\
 & =\bm{w}^{T}(\bm{w}'-\bm{w})/\left\Vert \bm{w}\right\Vert ^{2}\leq\left\Vert \bm{w}'-\bm{w}\right\Vert /\left\Vert \bm{w}\right\Vert ,\qquad\text{and}\\
\left\Vert A^{-T}ZA^{-1}-Z\right\Vert _{2} & \leq\left\Vert A^{-T}ZA^{-1}-A^{-T}Z\right\Vert _{2}+\left\Vert A^{-T}Z-Z\right\Vert _{2}\\
 & \leq\left\Vert A^{-T}\right\Vert _{2}\left\Vert Z\right\Vert _{2}\left\Vert A^{-1}-I\right\Vert _{2}+\left\Vert A^{-T}-I\right\Vert _{2}\left\Vert Z\right\Vert _{2}\\
 & =(1+\left\Vert A^{-1}\right\Vert _{2})\left\Vert A^{-1}-I\right\Vert _{2}\left\Vert Z\right\Vert _{2}.
\end{align*}
Combining the two parts, $\E[\Hess\widehat{f}_{\alpha}(\bm{w};\mathcal{D})]$
is Lipschitz on compact sets not containing $\bm{w}=\bm{0}$,
uniformly as $\alpha\downarrow0$. Note that for any dataset $\mathcal{D}$,
the function $\Hess\widehat{f}_{\alpha}(\bm{w};\mathcal{D})$
is Lipschitz in $\bm{w}$ with a Lipschitz constant that is
$\mathcal{O}(1/\alpha^{2})$ as $\alpha\downarrow0$.
\end{proof}

\subsubsection{Lipschitz constants for $\Hess\widehat{f}_{\alpha}(\bm{w};\mathcal{D})$} \label{subsec:Lips-const-Hessf}
\begin{lemma}
With probability one, if $C$ is a compact set not containing $\bm{0}$, we have
\[
\left\Vert \Hess\widehat{f}_{\alpha}(\bm{w};\mathcal{D})-\Hess\widehat{f}_{\alpha}(\bm{w}';\mathcal{D})\right\Vert _{2}/\left\Vert \bm{w}-\bm{w}'\right\Vert = \mathcal{O}\left(1+1/\sqrt{N\alpha^3}\right), \quad \forall\, \bm{w},\bm{w}' \in C,
\]
as $\alpha\downarrow0$.
\end{lemma}
\begin{proof}
We first need to estimate 
\begin{align*}
 & \Var\left[\Hess\widehat{f}_{\alpha}(\bm{w};\mathcal{D})-\Hess\widehat{f}_{\alpha}(\bm{w}';\mathcal{D})\right]\\
 & \quad{}=\Var\left[\frac{1}{N}\sum_{i=1}^{N}\frac{1}{\alpha}\left(\psi''\left(\frac{1-y_{i}\bm{x}_{i}^{T}\bm{w}}{\alpha}\right)-\psi''\left(\frac{1-y_{i}\bm{x}_{i}^{T}\bm{w}'}{\alpha}\right)\right)\bm{x}_{i}\bm{x}_{i}^{T}\right]\\
 & \quad{}=\frac{1}{(N\alpha)^{2}}\Var\left[\sum_{i=1}^{N}\left(\psi''\left(\frac{1-y_{i}\bm{x}_{i}^{T}\bm{w}}{\alpha}\right)-\psi''\left(\frac{1-y_{i}\bm{x}_{i}^{T}\bm{w}'}{\alpha}\right)\right)\bm{x}_{i}\bm{x}_{i}^{T}\right]\\
 & \quad{}=\frac{1}{N^{2}\alpha^{2}}\sum_{i=1}^{N}\Var\left[\left(\psi''\left(\frac{1-y_{i}\bm{x}_{i}^{T}\bm{w}}{\alpha}\right)-\psi''\left(\frac{1-y_{i}\bm{x}_{i}^{T}\bm{w}'}{\alpha}\right)\right)\bm{x}_{i}\bm{x}_{i}^{T}\right]\\
 & \quad{}=\frac{N}{N^{2}\alpha^{2}}\Var_{(\bm{X},Y)\sim p}\left[\left(\psi''\left(\frac{1-Y\bm{X}^{T}\bm{w}}{\alpha}\right)-\psi''\left(\frac{1-Y\bm{X}^{T}\bm{w}'}{\alpha}\right)\right)\bm{X}\bm{X}^{T}\right]\\
 & \quad{}\leq\frac{1}{N\alpha^{2}}\E_{(\bm{X},Y)\sim p}\left[\left(\psi''\left(\frac{1-Y\bm{X}^{T}\bm{w}}{\alpha}\right)-\psi''\left(\frac{1-Y\bm{X}^{T}\bm{w}'}{\alpha}\right)\right)^{2}\left\Vert \bm{X}\right\Vert ^{2}\right]\\
 & \quad{}=\frac{1}{N\alpha^{2}}\sum_{y\in\left\{ \pm1\right\} }\int_{\R^{m}}p(\bm{x},y)\left(\psi''\left(\frac{1-y\bm{x}^{T}\bm{w}}{\alpha}\right)-\psi''\left(\frac{1-y\bm{x}^{T}\bm{w}'}{\alpha}\right)\right)^{2}\left\Vert \bm{x}\right\Vert ^{2}d\bm{x}
\end{align*}
Now, noting that $\psi''(s)=\frac{1}{2}(1+s^{2})^{-3/2}$,
\[
\frac{\left(\psi''(s)-\psi''(t)\right)^{2}}{\psi''(s)^{2}+\psi''(t)^{2}}\leq C\,\left(s-t\right)^{2}\qquad\text{for all }s,t\in\R,
\]
with the value of $C=9/8=1.125$. Thus
\begin{align*}
 & \Var\left[\Hess\widehat{f}_{\alpha}(\bm{w};\mathcal{D})-\Hess\widehat{f}_{\alpha}(\bm{w}',\mathcal{D})\right]\\
 & \quad{}\leq\frac{1}{N\alpha^{2}}\sum_{y\in\left\{ \pm1\right\} }\int_{\R^{m}}p(\bm{x},y)\,C\left(\frac{y\bm{x}^{T}(\bm{w}-\bm{w}')}{\alpha}\right)^{2}\\
 & \qquad\qquad{}\times\left(\psi''\left(\frac{1-y\bm{x}^{T}\bm{w}}{\alpha}\right)^{2}+\psi''\left(\frac{1-y\bm{x}^{T}\bm{w}'}{\alpha}\right)^{2}\right)\left\Vert \bm{x}\right\Vert ^{2}d\bm{x}\\
 & \quad{}\leq\frac{C\left\Vert \bm{w}-\bm{w}'\right\Vert ^{2}}{N\alpha^{4}}\sum_{y\in\left\{ \pm1\right\} }\int_{\R^{m}}p(\bm{x},y)\,\\
 & \qquad\qquad{}\times\left(\psi''\left(\frac{1-y\bm{x}^{T}\bm{w}}{\alpha}\right)^{2}+\psi''\left(\frac{1-y\bm{x}^{T}\bm{w}'}{\alpha}\right)^{2}\right)\left\Vert \bm{x}\right\Vert ^{4}d\bm{x}.
\end{align*}
Since $\bm{x}\mapsto p(\bm{x},y)$ is continuous,
bounded, with bounded support and $u\mapsto\psi''(u)^{2}$ is integrable,
the integrals
\[
\int_{\R^{m}}p(\bm{x},y)\,\left(\psi''\left(\frac{1-y\bm{x}^{T}\bm{w}}{\alpha}\right)^{2}+\psi''\left(\frac{1-y\bm{x}^{T}\bm{w}'}{\alpha}\right)^{2}\right)\left\Vert \bm{x}\right\Vert ^{4}d\bm{x}=\mathcal{O}(\alpha),
\]
with a hidden constant that depends only on $p$ and $\psi$. Therefore,
\begin{align*}
\Var\left[\Hess\widehat{f}_{\alpha}(\bm{w};\mathcal{D})-\Hess\widehat{f}_{\alpha}(\bm{w}',\mathcal{D})\right] & =\frac{\left\Vert \bm{w}-\bm{w}'\right\Vert ^{2}}{N\alpha^{4}}\mathcal{O}(\alpha)=\mathcal{O}\left(\frac{1}{N\alpha^{3}}\right)\left\Vert \bm{w}-\bm{w}'\right\Vert ^{2}.
\end{align*}
Noting that
\begin{align*}
 & \E\left[\left\Vert \Hess\widehat{f}_{\alpha}(\bm{w};\mathcal{D})-\Hess\widehat{f}_{\alpha}(\bm{w}',\mathcal{D})\right\Vert _{F}^{2}\right]\\
 & \quad{}=\left\Vert \E\left[\Hess\widehat{f}_{\alpha}(\bm{w};\mathcal{D})-\Hess\widehat{f}_{\alpha}(\bm{w}',\mathcal{D})\right]\right\Vert _{F}^{2}+\Var\left[\Hess\widehat{f}_{\alpha}(\bm{w};\mathcal{D})-\Hess\widehat{f}_{\alpha}(\bm{w}',\mathcal{D})\right]\\
 & \quad{}=\mathcal{O}(1)\left\Vert \bm{w}-\bm{w}'\right\Vert ^{2}+\mathcal{O}\left(\frac{1}{N\alpha^{3}}\right)\left\Vert \bm{w}-\bm{w}'\right\Vert ^{2},
\end{align*}
we can then use Jensen's inequality to obtain
\[
\E\left[\left\Vert \Hess\widehat{f}_{\alpha}(\bm{w};\mathcal{D})-\Hess\widehat{f}_{\alpha}(\bm{w}',\mathcal{D})\right\Vert _{F}\right]=\mathcal{O}\left(1+\frac{1}{\left(N\alpha\right)^{1/2}\,\alpha}\right)\left\Vert \bm{w}-\bm{w}'\right\Vert .
\]
By Chebyshev's bound, 
\[
\left\Vert \Hess\widehat{f}_{\alpha}(\bm{w};\mathcal{D})-\Hess\widehat{f}_{\alpha}(\bm{w}',\mathcal{D})\right\Vert _{F}=\mathcal{O}\left(1+\frac{1}{\left(N\alpha\right)^{1/2}\,\alpha}\right)\left\Vert \bm{w}-\bm{w}'\right\Vert 
\]
with probability one. 
\end{proof}

\subsubsection{Bound $0\prec aI\preceq\left\Vert \bm{w}\right\Vert \protect\E\left[\Hess\widehat{f}_{\alpha}(\bm{w};\mathcal{D})\right]\preceq bI$ for all sufficiently small $\alpha>0$} \label{subsec:Uniform-bdd-condEHess}
From the assumptions that the probability density of the data points $(\bm{x}_{i},y_{i})$ is Lipschitz continuous and the density has bounded support, we can show that there is a bound on the 2-norm condition number of $\E[\Hess\widehat{f}_{\alpha}(\bm{w};\mathcal{D})]$ provided $\bm{w}$ is in a bounded region and bounded away from zero. This bound, however, can grow exponentially in the dimension for a fixed Lipschitz constant.
\begin{lemma}
Let  $W := \big\{ \bm{w}\in\R^{m} \colon \sum_{y\in\left\{ \pm1\right\} }\int_{\left\{ \bm{x}\mid y=\bm{x}^{T}\bm{w}\right\} }p(\bm{x},y)\,dS(\bm{x})>0 \big\}$ and let $W_{0}$ be a compact subset of $W$. Then, there are constants $b > a > 0$ such that
\[
0\prec aI\preceq\left\Vert \bm{w}\right\Vert \E\left[\Hess\widehat{f}_{\alpha}(\bm{w};\mathcal{D})\right]\preceq bI, \qquad \forall \bm{w} \in W_0,
\]
for all sufficiently small $\alpha>0$.
\end{lemma}
\begin{proof}
For any $\bm{z}$, 
\begin{align}
 & \bm{z}^{T}\E[\Hess\widehat{f}_{\alpha}(\bm{w};\mathcal{D})]\bm{z}\nonumber \\
 & \qquad{}=\bm{z}^{T}\sum_{y\in\left\{ \pm1\right\} }\int_{\R^{m}}p(\bm{x},y)\frac{1}{\alpha}\psi''\left(\frac{1-y\bm{x}^{T}\bm{w}}{\alpha}\right)\bm{x}\bm{x}^{T}\,d\bm{x}\,\bm{z}+\lambda\bm{z}^{T}\bm{z}\nonumber \\
 & \qquad{}=\sum_{y\in\left\{ \pm1\right\} }\int_{\R^{m}}p(\bm{x},y)\frac{1}{\alpha}\psi''\left(\frac{1-y\bm{x}^{T}\bm{w}}{\alpha}\right)(\bm{z}^{T}\bm{x})^{2}\,d\bm{x}+\lambda\bm{z}^{T}\bm{z}\nonumber \\
 & \qquad{}\to \frac{1}{\|\bm{w}\|}\sum_{y\in\left\{ \pm1\right\} }\int_{\left\{ \bm{x}\mid y=\bm{x}^{T}\bm{w}\right\} }p(\bm{x},y)(\bm{z}^{T}\bm{x})^{2}\,dS(\bm{x}) +\lambda\bm{z}^{T}\bm{z}\qquad\text{as }\alpha\downarrow0.\label{eq:lim-z.Hess.z}
\end{align}
Let 
\[
h(\bm{z},\bm{w})=\sum_{y\in\left\{ \pm1\right\} }\int_{\left\{ \bm{x}\mid y=\bm{x}^{T}\bm{w}\right\} }p(\bm{x},y)(\bm{z}^{T}\bm{x})^{2}\,dS(\bm{x}).
\]
By continuity of $h(\bm{z},\bm{w})$ for $\bm{w}\neq\bm{0}$, for any $0<r_{0}<R_{0}$ we can set 
\[
b_{0}=2\,\max\big\{ h(\bm{z},\bm{w})\mid\left\Vert \bm{z}\right\Vert =1,\,r_{0}\leq\left\Vert \bm{w}\right\Vert \leq R_{0}\big\}.
\]
Note that if $\lambda>0$ we can take $a_{0}=0$ as $h(\bm{z},\bm{w})\geq0$ for any $\bm{z}$ and $\bm{w}$. Under the assumptions of Section~\ref{subsec:Assumptions}, either $\lambda>0$ or the optimal $\bm{w}$ implies $\sum_{y\in\left\{ \pm1\right\} }\int_{\left\{ \bm{x}\mid y=\bm{x}^{T}\bm{w}\right\} }p(\bm{x},y)\,dS(\bm{x})>0$.
We now focus on the case where $\lambda=0$. Let $W$ and $W_{0}$ be sets as described in the statement of the Lemma; by continuity of $p(\bm{x},y)$, $\bm{w}\in W$ is equivalent to there being a point $(\bm{x},y)$ where $\bm{x}^{T}\bm{w}=y=\pm1$ and $p(\bm{x},y)>0$. By continuity of $p(\bm{x},y)$, $W$ is an open set. Note that $\bm{0}\not\in W$ and $\bm{w}^{*}\in W$. It follows that $\min\left\{ h(\bm{z},\bm{w})\mid\left\Vert \bm{z}\right\Vert =1,\,\bm{w}\in W_{0}\right\} >0$.
Set 
\[
a_{0}=\frac{1}{2}\min\left\{ \,h(\bm{z},\bm{w})\mid\left\Vert \bm{z}\right\Vert =1,\,\bm{w}\in W_{0}\,\right\} >0.
\]
Since $\bm{0}\not\in W$ as $\left\{ \bm{x}\mid y=\bm{x}^{T}\bm{w}\right\} =\emptyset$ for $\bm{w}=\bm{0}$ and $y=\pm1$, we have $\bm{0}\not\in W_{0}$. Thus, we can take $r_{0}=\min_{\bm{w}\in W_{0}}\left\Vert \bm{w}\right\Vert$ and $R_{0}=\max_{\bm{w}\in W_{0}}\left\Vert \bm{w}\right\Vert$. We have $h(\bm{z},\bm{w})\leq\frac{1}{2}b_{0}\left\Vert \bm{z}\right\Vert ^{2}$ for all $\bm{w}\in W_{0}$. Hence, for any $\bm{w}\in W_{0}$, we get
\begin{align*}
\lim_{\alpha\downarrow0}\lambda_{\min}\left(\E\big[\Hess\widehat{f}_{\alpha}(\bm{w};\mathcal{D})\big] \right)
&\geq2a_{0}+\lambda\left\Vert \bm{w}\right\Vert, \\
\lim_{\alpha\downarrow0}\lambda_{\max}\left(\E\big[\Hess\widehat{f}_{\alpha}(\bm{w};\mathcal{D})\big]\right)
&\leq\frac{1}{2}b_{0}+\lambda\left\Vert \bm{w}\right\Vert.
\end{align*}
Therefore, for sufficiently small $\alpha>0$, we have
\[
(a_{0}+\lambda\left\Vert \bm{w}\right\Vert )I\preceq\left\Vert \bm{w}\right\Vert \E\big[\Hess\widehat{f}_{\alpha}(\bm{w};\mathcal{D})\big]\preceq(b_{0}+\lambda\left\Vert \bm{w}\right\Vert )I, \quad \forall \bm{w}\in W_{0}.
\]
Setting $a=a_{0}+\lambda r_{0}>0$ and $b=b_{0}+\lambda R_{0}$ gives
the bounds we seek.%
\end{proof}
We can choose $W_{0}$ to be a compact subset of $W$ containing $\bm{w}^{*}$
in its interior. Then any sufficiently close $\bm{w}\approx\bm{w}^{*}$
is also in $W_{0}$. We can take, for example, $W_{0}=\left\{ \,\bm{w}\mid d(\bm{w},\R^{m}\backslash W)\geq\delta\,\right\} $
for some sufficiently small $\delta>0$ such as $\delta=\frac{1}{2}d(\bm{w}^{*},\R^{m}\backslash W)$.

\subsubsection{Variation of $\Hess\widehat{f}_{\alpha}(\bm{w};\mathcal{D})$} \label{subsec:Var-Hess}
While we can represent the behavior of $\E[\Hess\widehat{f}_{\alpha}(\bm{w};\mathcal{D})]$ in terms of integrals that behave nicely, the samples in the dataset can exhibit randomness that results in significant variation that does not appear in the above bounds (except that $\Hess\widehat{f}_{\alpha}(\bm{w};\mathcal{D})$ is Lipschitz with constant $\mathcal{O}(1/\alpha^{2})$). This additional variation can be important even for $N\alpha\gg1$ and $\alpha\ll1$. With $\alpha$ small and $\left\Vert \bm{w}-\bm{w}'\right\Vert _{2}\gg\alpha$ the effective supports of the function $\psi''((1-y\bm{w}^{T}\bm{x})/\alpha)$ can be disjoint or have small overlap meaning that the data points used for computing $\Hess\widehat{f}_{\alpha}(\bm{w};\mathcal{D})$ and $\Hess\widehat{f}_{\alpha}(\bm{w}';\mathcal{D})$ may be nearly if not entirely disjoint. In these cases, the main connection between $\Hess\widehat{f}_{\alpha}(\bm{w};\mathcal{D})$ and $\Hess\widehat{f}_{\alpha}(\bm{w}',\mathcal{D})$ is through the probability distribution of the data points rather than the specific data points themselves. The variance of $\Hess\widehat{f}_{\alpha}(\bm{w};\mathcal{D})$ can give good estimates for this variation.

To be more precise, we let $\Var[Z]=\E\left[(Z-\E[Z])\bullet(Z-\E[Z])\right]$ for a matrix-value random variable, where $A\bullet B=\trace(A^{T}B)=\sum_{i,j}a_{ij}b_{ij}$ is the Frobenius inner product \cite[p.~332, Ex.~24]{Hor85}. Note that with this definition, $\Var[Z]=\E\left[\left\Vert Z-\E[Z]\right\Vert _{F}^{2}\right]$ where $\left\Vert A\right\Vert _{F}=\sqrt{A\bullet A}=\sqrt{\sum_{ij}\left|a_{ij}\right|^{2}}$ is the Frobenius matrix norm. By convexity of the function $u\mapsto u^{2}$ and Jensen's inequality \cite[pp.~44-45]{Lie2001}, we have
\[
\E\left[\left\Vert Z-\E[Z]\right\Vert _{F}\right]^{2}\leq\E\left[\left\Vert Z-\E[Z]\right\Vert _{F}^{2}\right]=\Var[Z].
\]
\begin{lemma}
$\Var\left[\Hess\widehat{f}_{\alpha}(\bm{w};\mathcal{D})\right]=\mathcal{O}(1/(N\alpha))$
as $\alpha\downarrow0$.
\end{lemma}

\begin{proof}
Note that
\begin{align*}
\Var\left[\Hess\widehat{f}_{\alpha}(\bm{w};\mathcal{D})\right] & =\frac{1}{N^{2}\alpha^{2}}\sum_{i=1}^{N}\Var_{(\bm{X},Y)\sim p}\left[\psi''\left(\frac{1-Y\bm{w}^{T}\bm{X}}{\alpha}\right)\bm{X}\bm{X}^{T}\right]\\
 & =\frac{1}{N\alpha^{2}}\Var_{(\bm{X},Y)\sim p}\left[\psi''\left(\frac{1-Y\bm{w}^{T}\bm{X}}{\alpha}\right)\bm{X}\bm{X}^{T}\right].
\end{align*}
Now, 
\begin{align*}
 & \left\Vert \E_{(\bm{X},Y)\sim p}\left[\psi''\left(\frac{1-Y\bm{w}^{T}\bm{X}}{\alpha}\right)\bm{X}\bm{X}^{T}\right]\right\Vert _{F}\\
 & \qquad{}\leq\sum_{y\in\left\{ \pm1\right\} }\int_{\R^{m}}p(\bm{x},y)\,\psi''\left(\frac{1-y\bm{w}^{T}\bm{x}}{\alpha}\right)\left\Vert \bm{x}\bm{x}^{T}\right\Vert _{F}\,d\bm{x}\\
 & \qquad{}=\sum_{y\in\left\{ \pm1\right\} }\int_{-\infty}^{+\infty}\int_{\left\{ \bm{w}\right\} ^{\perp}}p(s\widehat{\bm{w}}+\bm{v})\,\psi''\left(\frac{1-ys\left\Vert \bm{w}\right\Vert _{2}}{\alpha}\right)\left\Vert s\widehat{\bm{w}}+\bm{v}\right\Vert _{2}^{2}d\bm{v}\,ds\\
 & \qquad\qquad\qquad(\bm{x}=s\widehat{\bm{w}}+\bm{v}\text{ with }\bm{v}\perp\bm{w})\\
 & \qquad{}\leq\sum_{y\in\left\{ \pm1\right\} }\int_{-\infty}^{+\infty}p_{\bm{w}}(s)\,\psi''\left(\frac{1-ys\left\Vert \bm{w}\right\Vert _{2}}{\alpha}\right)\,R^{2}\,ds=\mathcal{O}(\alpha).
\end{align*}
On the other hand,
\begin{align*}
 & \E_{(\bm{X},Y)\sim p}\left[\psi''\left((1-Y\bm{w}^{T}\bm{X})/\alpha\right)\bm{X}\bm{X}^{T}\bullet\psi''\left((1-Y\bm{w}^{T}\bm{X})/\alpha\right)\bm{X}\bm{X}^{T}\right]\\
 & \qquad{}=\E_{(\bm{X},Y)\sim p}\left[\psi''\left((1-Y\bm{w}^{T}\bm{X})/\alpha\right)^{2}\left\Vert \bm{X}\right\Vert _{2}^{4}\right]\\
 & \qquad{}\leq R^{4}\,\E_{(\bm{X},Y)\sim p}\left[\psi''\left((1-Y\bm{w}^{T}\bm{X})/\alpha\right)^{2}\right]=\mathcal{O}(\alpha)
\end{align*}
using similar integration techniques to the previous computation.
Thus,
\[
Z =\psi''\left((1-Y\bm{w}^{T}\bm{X})/\alpha\right)\bm{X}\bm{X}^{T}
\]
implies 
\[
\Var\left[Z\right] =\E\left[Z\bullet Z\right]-\E\left[Z\right]\bullet\E\left[Z\right]\leq\E\left[Z\bullet Z\right]=\mathcal{O}(\alpha).
\]
Note that disregarding the $\E\left[Z\right]\bullet\E\left[Z\right]$
term does not come at a great cost since this term is $\mathcal{O}(\alpha)\times\mathcal{O}(\alpha)=\mathcal{O}(\alpha^{2})$
and we are interested in $\alpha\ll1$. Thus,
\[
\Var\left[\Hess\widehat{f}_{\alpha}(\bm{w};\mathcal{D})\right]=\frac{1}{N\alpha^{2}}\mathcal{O}(\alpha)=\mathcal{O}\left(\frac{1}{N\alpha}\right)
\]
as desired.
\end{proof}
An easy consequence by Jensen's inequality is that
\begin{align}
\E\left[\Big\Vert \Hess\widehat{f}_{\alpha}(\bm{w};\mathcal{D})-\E\big[\Hess\widehat{f}_{\alpha}(\bm{w};\mathcal{D})\big]\Big\Vert _{F}\right] &\leq\sqrt{\Var\big[\Hess\widehat{f}_{\alpha}(\bm{w};\mathcal{D})\big]} \nonumber\\
&=\mathcal{O}\left(\frac{1}{\sqrt{N\alpha}}\right),\label{eq:Var-Hessf}
\end{align}
this is asymptotically small, provided $N\alpha\gg1$.

\subsubsection{Variance of gradients} \label{subsec:Variance-of-gradients}
Another important aspect is the estimation of the variance of the gradients $\Var\big[\nabla\widehat{f}_{\alpha}(\bm{w};\mathcal{D})\big]$.
In particular, since $\Var\left[\bm{Z}\right]$ is the variance\textendash covariance
matrix, $\trace\,\Var\left[\bm{Z}\right]=\E\big[\left\Vert \bm{Z}-\E\left[\bm{Z}\right]\right\Vert _{2}^{2}\big]$
gives an estimate of the variation of the gradient across different
samples $\mathcal{D}$. 
\begin{lemma}
$\trace \Var\left[\nabla\widehat{f}_{\alpha}(\bm{w};\mathcal{D})\right]=\mathcal{O}(1/N)$
as $\alpha\downarrow0$.
\end{lemma}
\begin{proof}
To bound $\trace\,\Var\left[\nabla\widehat{f}_{\alpha}(\bm{w};\mathcal{D})\right]$,
note that under the assumption that each data point $(\bm{x}_{i},y_{i})\in\R^{m}\times\left\{ \pm1\right\} $
is chosen randomly and independently according a probability distribution
with probability density function $p$,
\begin{align*}
\Var\left[\nabla\widehat{f}_{\alpha}(\bm{w};\mathcal{D})\right] & =\Var\left[N^{-1}\sum_{i=1}^{N}\psi'\left(\frac{1-\bm{w}^{T}\bm{x}_{i}y_{i}}{\alpha}\right)\bm{x}_{i}y_{i}\right]\\
 & =N^{-2}\sum_{i=1}^{N}\Var\left[\psi'\left(\frac{1-\bm{w}^{T}\bm{x}_{i}y_{i}}{\alpha}\right)\bm{x}_{i}y_{i}\right]\\
 & =N^{-1}\Var_{(\bm{X},Y)\sim p}\left[\psi'\left(\frac{1-\bm{w}^{T}\bm{X}Y}{\alpha}\right)\bm{X}Y\right].
\end{align*}
Since $\psi'$ is a bounded function: $\left|\psi'(s)\right|\leq1$
for all $s$, we have
\begin{equation}
\left\Vert \Var\left[\nabla\widehat{f}_{\alpha}(\bm{w};\mathcal{D})\right]\right\Vert _{2}\leq N^{-1}\,\E_{(\bm{X},Y)\sim p}\left[\left\Vert \bm{X}\right\Vert _{2}^{2}\right]=\mathcal{O}(1/N),\label{eq:Var-gradients}
\end{equation}
as we wanted.
\end{proof}

\subsubsection{Effects of changing $\alpha$} \label{subsec:Changing-alpha}
Changing $\alpha$ via $\alpha\gets\beta\alpha$ with $0<\beta<1$
results in a different smooth optimization problem. If we nearly have
a solution to the problem with smoothing parameter $\alpha$, it should
also not be far from the solution for smoothing parameter $\beta\,\alpha$
for $0<\beta<1$ fixed.
\begin{lemma}
If $\left\Vert \nabla\widehat{f}_{\alpha}(\bm{w};\mathcal{D})\right\Vert _{2}=\mathcal{O}(\alpha)$,
and $1\gg\alpha\gg N^{-1}$ then with probability one \emph{$\left\Vert \nabla\widehat{f}_{\alpha\beta}(\bm{w};\mathcal{D})\right\Vert _{2}=\mathcal{O}(\alpha)$.}
\end{lemma}

\begin{proof}
To see why, we compute
\begin{align*}
 & \E[\nabla\widehat{f}_{\alpha}(\bm{w};\mathcal{D})-\nabla\widehat{f}_{\beta\alpha}(\bm{w};\mathcal{D})]=\\
 & \qquad\sum_{y\in\left\{ \pm1\right\} }\int_{\R^{m}}p(\bm{x},y)\,\left[\psi'\left(\frac{1-y\bm{w}^{T}\bm{x}}{\alpha}\right)-\psi'\left(\frac{1-y\bm{w}^{T}\bm{x}}{\beta\alpha}\right)\right]y\bm{x}\,d\bm{x}.
\end{align*}
Writing $\bm{x}=s\widehat{\bm{w}}+\bm{v}$
with $\bm{v}\perp\bm{w}$ and $\widehat{\bm{w}}=\bm{w}/\left\Vert \bm{w}\right\Vert _{2}$,
the integral for \\
$\left\Vert \E[\nabla\widehat{f}_{\alpha}(\bm{w};\mathcal{D})-\nabla\widehat{f}_{\beta\alpha}(\bm{w};\mathcal{D})]\right\Vert _{2}$
becomes
\begin{align*}
 & \left\Vert \sum_{y\in\left\{ \pm1\right\} }\int_{-\infty}^{+\infty}\int_{\left\{ \bm{w}\right\} ^{\perp}}p(s\widehat{\bm{w}}+\bm{v},y)\,\left[\psi'\left(\frac{1-ys\left\Vert \bm{w}\right\Vert _{2}}{\alpha}\right)-\psi'\left(\frac{1-ys\left\Vert \bm{w}\right\Vert _{2}}{\beta\alpha}\right)\right]\times{}\right.\\
 & \qquad\qquad\qquad{}\times\left.\vphantom{\left[\left(\frac{1-ys\left\Vert \bm{w}\right\Vert _{2}}{\alpha}\right)\right]}y(s\widehat{\bm{w}}+\bm{v})\,d\bm{v}\,ds\right\Vert _{2}\\
 & \qquad{}\leq\sum_{y\in\left\{ \pm1\right\} }\int_{-\infty}^{+\infty}p_{\bm{w}}(s,y)\left|\psi'\left(\frac{1-ys\left\Vert \bm{w}\right\Vert _{2}}{\alpha}\right)-\psi'\left(\frac{1-ys\left\Vert \bm{w}\right\Vert _{2}}{\beta\alpha}\right)\right|R\,ds\\
 & \qquad{}\leq\sum_{y\in\left\{ \pm1\right\} }\int_{-\infty}^{+\infty}p_{\bm{w}}\left(\frac{1-\alpha t}{y\left\Vert \bm{w}\right\Vert _{2}},y\right)\left|\psi'\left(t\right)-\psi'\left(t/\beta\right)\right|R\,\frac{\alpha}{\left\Vert \bm{w}\right\Vert _{2}}dt\\
 & \qquad\qquad\qquad(t=(1-ys\left\Vert \bm{w}\right\Vert _{2})/\alpha)\\
 & \qquad{}\leq\sum_{y\in\left\{ \pm1\right\} }\max_{s}p_{\bm{w}}(s,y)\,\frac{R\alpha}{\left\Vert \bm{w}\right\Vert _{2}}\int_{-\infty}^{+\infty}\left|\psi'\left(t\right)-\psi'\left(t/\beta\right)\right|\,dt=\mathcal{O}(\alpha)
\end{align*}
as $\left|\psi'\left(t\right)-\psi'\left(t/\beta\right)\right|=\mathcal{O}(1/t^{2})$
as $t\to\pm\infty$ for $\psi(t)=\frac{1}{2}\left(t+\sqrt{1+t^{2}}\right)$.
Note, that other suitable functions like $\psi(t)=\ln(1+\exp(t))$
have even faster decay of $\left|\psi'\left(t\right)-\psi'\left(t/\beta\right)\right|$
as $t\to\pm\infty$, so that these arguments still apply.

We can obtain bounds on $\E\left[\left\Vert \nabla f_{\alpha}(\bm{w};\mathcal{D})-\nabla f_{\beta\alpha}(\bm{w};\mathcal{D})\right\Vert _{2}\right]$
as follows: from Jensen's inequality, for any random variable $\bm{Z}$
with values in $\R^{m}$ and finite variance we have $\E\left[\left\Vert \bm{Z}-\E\left[\bm{Z}\right]\right\Vert \right]^{2}\leq\E\left[\left\Vert \bm{Z}-\E\left[\bm{Z}\right]\right\Vert ^{2}\right]=\trace\Var\left[\bm{Z}\right]$.
Now
\begin{align*}
 & \trace\Var\left[\nabla f_{\alpha}(\bm{w};\mathcal{D})-\nabla f_{\beta\alpha}(\bm{w};\mathcal{D})\right]\\
 & \qquad{}=\trace\Var\left[N^{-1}\sum_{i=1}^{N}\left(\psi'\left(\frac{1-\bm{w}^{T}\bm{x}_{i}y_{i}}{\alpha}\right)\bm{x}_{i}y_{i}-\psi'\left(\frac{1-\bm{w}^{T}\bm{x}_{i}y_{i}}{\alpha\beta}\right)\bm{x}_{i}y_{i}\right)\right]\\
 & \qquad{}=N^{-2}\sum_{i=1}^{N}\trace\Var\left[\left\{ \psi'\left(\frac{1-\bm{w}^{T}\bm{x}_{i}y_{i}}{\alpha}\right)-\psi'\left(\frac{1-\bm{w}^{T}\bm{x}_{i}y_{i}}{\alpha\beta}\right)\right\} \bm{x}_{i}y_{i}\right]\\
 & \qquad{}=N^{-1}\trace\Var_{(\bm{X},Y)\sim p}\left[\left\{ \psi'\left(\frac{1-\bm{w}^{T}\bm{X}Y}{\alpha}\right)-\psi'\left(\frac{1-\bm{w}^{T}\bm{X}Y}{\alpha\beta}\right)\right\} \bm{X}Y\right]\\
 & \qquad{}\leq N^{-1}\sum_{y\in\left\{ \pm1\right\} }\int_{\R^{m}}p(\bm{x},y)\,\left\{ \psi'\left(\frac{1-\bm{w}^{T}\bm{x}y}{\alpha}\right)-\psi'\left(\frac{1-\bm{w}^{T}\bm{x}y}{\alpha\beta}\right)\right\} ^{2}\left\Vert \bm{x}\right\Vert _{2}^{2}\,d\bm{x}.
\end{align*}
Writing $\bm{x}=s\widehat{\bm{w}}+\bm{v}$
with $\bm{v}\perp\bm{w}$ and $\widehat{\bm{w}}=\bm{w}/\left\Vert \bm{w}\right\Vert _{2}$,
we get
\begin{align*}
 & \trace\Var\left[\nabla f_{\alpha}(\bm{w};\mathcal{D})-\nabla f_{\beta\alpha}(\bm{w};\mathcal{D})\right]\\
 & \qquad{}\leq N^{-1}\sum_{y\in\left\{ \pm1\right\} }\int_{-\infty}^{+\infty}\int_{\left\{ \bm{w}\right\} ^{\perp}}p(s\widehat{\bm{w}}+\bm{v})\\
 & \qquad\qquad{}\left\{ \psi'\left(\frac{1-\left\Vert \bm{w}\right\Vert _{2}sy}{\alpha}\right)-\psi'\left(\frac{1-\left\Vert \bm{w}\right\Vert _{2}sy}{\alpha\beta}\right)\right\} ^{2}\left(s^{2}+\left\Vert \bm{v}\right\Vert _{2}^{2}\right)\,d\bm{v}\,ds.
\end{align*}
Let $\theta(s,y,\bm{w})=\int_{\left\{ \bm{w}\right\} ^{\perp}}\left(s^{2}+\left\Vert \bm{v}\right\Vert _{2}^{2}\right)p(s\widehat{\bm{w}}+\bm{v})\,d\bm{v}$,
which under our assumptions is a bounded continuous function for $\bm{w}\neq0$,
and is zero for $\left|s\right|\geq R$. Then
\begin{align*}
 & \trace\Var\left[\nabla\widehat{f}_{\alpha}(\bm{w};\mathcal{D})-\nabla\widehat{f}_{\beta\alpha}(\bm{w};\mathcal{D})\right]\\
 & \qquad{}\leq N^{-1}\int_{-\infty}^{+\infty}\theta(s,y,\bm{w})\left\{ \psi'\left(\frac{1-\left\Vert \bm{w}\right\Vert _{2}sy}{\alpha}\right)-\psi'\left(\frac{1-\left\Vert \bm{w}\right\Vert _{2}sy}{\alpha\beta}\right)\right\} ^{2}\,ds\\
 & \qquad{}\leq N^{-1}\theta_{\max}\int_{-\infty}^{+\infty}\left\{ \psi'\left(\frac{1-\left\Vert \bm{w}\right\Vert _{2}sy}{\alpha}\right)-\psi'\left(\frac{1-\left\Vert \bm{w}\right\Vert _{2}sy}{\alpha\beta}\right)\right\} ^{2}\,ds.
\end{align*}
Setting $\widetilde{s}=(1-\left\Vert \bm{w}\right\Vert _{2}sy)/\alpha$
we see that after changing variables we have
\begin{align*}
 & \trace\Var\left[\nabla\widehat{f}_{\alpha}(\bm{w};\mathcal{D})-\nabla\widehat{f}_{\beta\alpha}(\bm{w};\mathcal{D})\right]\\
 & \qquad{}\leq\alpha N^{-1}\theta_{\max}\int_{-\infty}^{+\infty}\left\{ \psi'\left(\widetilde{s}\right)-\psi'\left(\widetilde{s}/\beta\right)\right\} ^{2}\,d\widetilde{s}.
\end{align*}
Now $\E\left[\left\Vert \bm{Z}\right\Vert _{2}\right]\leq\E\left[\left\Vert \bm{Z}-\E\left[\bm{Z}\right]\right\Vert _{2}+\left\Vert \E\left[\bm{Z}\right]\right\Vert _{2}\right]\leq\trace\Var\left[\bm{Z}\right]^{1/2}+\left\Vert \E\left[\bm{Z}\right]\right\Vert _{2}$.
Applied to this situation we have 
\begin{align*}
\E\left[\left\Vert \nabla\widehat{f}_{\alpha}(\bm{w};\mathcal{D})-\nabla\widehat{f}_{\beta\alpha}(\bm{w};\mathcal{D})\right\Vert _{2}\right] & \leq\mathcal{O}(\alpha^{1/2}N^{-1/2})+\mathcal{O}(\alpha)\\
 & =\mathcal{O}(\alpha(1+1/(\alpha N)^{1/2})).
\end{align*}
In the regime $N^{-1}\leq\alpha\ll1$, this is just becomes
\begin{align*}
\E\left[\left\Vert \nabla\widehat{f}_{\alpha}(\bm{w};\mathcal{D})-\nabla\widehat{f}_{\beta\alpha}(\bm{w};\mathcal{D})\right\Vert _{2}\right] & =\mathcal{O}(\alpha),\\
\trace\Var\left[\nabla\widehat{f}_{\alpha}(\bm{w};\mathcal{D})-\nabla\widehat{f}_{\beta\alpha}(\bm{w};\mathcal{D})\right] & =\mathcal{O}(\alpha/N).
\end{align*}
If $\alpha\gg N^{-1}$ then this implies that the standard deviation
of the difference is less than the size of the expected value of the
difference. So if $\left\Vert \nabla\widehat{f}_{\alpha}(\bm{w};\mathcal{D})\right\Vert _{2}=\mathcal{O}(\alpha)$,
we also have \emph{$\left\Vert \nabla\widehat{f}_{\alpha\beta}(\bm{w};\mathcal{D})\right\Vert _{2}=\mathcal{O}(\alpha)+\mathcal{O}(\alpha/N)^{1/2}=\mathcal{O}(\alpha)+\mathcal{O}(\alpha^{2})^{1/2}=\mathcal{O}(\alpha)$
}provided $\alpha\beta\geq N^{-1}$. 
\end{proof}

\subsubsection{Convergence of Newton's method in the midgame} \label{subsec:midgame-convce-Newton}
We wish to show that with probability one, the guarded or damped Newton method used here only requires $\mathcal{O}(1)$ iterations each time we reduce $\alpha$: $\alpha\gets\beta\alpha$ (Algorithm~\ref{alg:Adjust-active-set}, line~5). The guarded or damped Newton method has been analyzed in, for example, Nesterov and Nemirovskii \cite{Nes94} and Boyd, Boyd, and Vandenberghe \cite[Chap.~9]{Boy2004}. The standard Newton method has also been analyzed in Kantorovich \cite{Kan49} in which precise conditions are given for the well-known quadratic convergence of this method. 

Boyd, Boyd, and Vandenberghe \cite[Chap.~9]{Boy2004} give a bound on the number of guarded Newton iterations in terms of $\widehat{f}_{\alpha}(\bm{w}_{0},\mathcal{D})-\inf_{\bm{w}}\widehat{f}_{\alpha}(\bm{w};\mathcal{D})$ and the the target accuracy $\epsilon>0$, using bounds $aI\preceq\text{Hess}_{\bm{w}}\,\widehat{f}_{\alpha}(\bm{w};\mathcal{D})\preceq bI$ with positive constants $a$ and $b$, and a Lipschitz constant $L$ for $\text{Hess}_{\bm{w}}\,\widehat{f}_{\alpha}(\bm{w};\mathcal{D})$. Nesterov and Nemirovskii \cite{Nes94} instead use self-concordance (there is a constant $C$ such that for any $\bm{w}$ and $\bm{d}$, the function $k(t):=\widehat f_\alpha(\bm{w}+t\bm{d};\mathcal{D})$ has $|k'''(t)|/k''(t)^{3/2}\leq C$ for all $t$) to obtain bounds on the number of guarded Newton iterations needed. However, for our case, the function $\psi(u)=\frac{1}{2}\big[u+\sqrt{1+u^{2}}\big]$ is not self-concordant as $\left|\psi'''(u)\right|/\psi''(u)=3\sqrt{2}\left|u\right|(1+u^{2})^{-1/4}\to\infty$
as $u\to\pm\infty$. 
\begin{lemma}
If $1\gg\alpha\gg N^{-1}$ and $0<\eta<1$ is fixed, then with probability one, $\mathcal{O}(1+\log\log(1/\alpha))$ for $1 \gg \alpha$ guarded Newton steps are sufficient to satisfy the stopping criterion of Algorithm~\ref{alg:smoothed-l1-Newton-step}, line~6. If, in addition, $\eta^{2}\,\alpha\gg N^{-1}$, then only one Newton step
without any further line search is needed.
\end{lemma}

\begin{proof}
Boyd, Boyd, and Vandenberghe \cite[p.~491]{Boy2004} show that guarded
Newton's method with Armijo backtracking to find $\bm{w}$
where $\left|\bm{d}^{T}\nabla\widehat{f}_{\alpha}(\bm{w};\mathcal{D})\right|<\epsilon$
where $\bm{d}$ is the Newton step from $\bm{w}$,
uses no more than $\mathcal{O}(1+\gamma(\widehat{f}_{\alpha}(\bm{w};\mathcal{D})-\inf_{\bm{w}'}\widehat{f}_{\alpha}(\bm{w}',\mathcal{D})))+\mathcal{O}(\log\log(\epsilon_{0}/\epsilon))$
Newton steps. Here $\gamma=(b')^{2}L^{2}/(a')^{5}$ where $0\prec a'I\preceq\Hess\widehat{f}_{\alpha}(\bm{w};\mathcal{D})\preceq b'I$,
$L$ is the Lipschitz constant for $\bm{w}\mapsto\Hess\widehat{f}_{\alpha}(\bm{w};\mathcal{D})$,
and $\epsilon_{0}=(a')^{3}/L^{2}$. From Sections~\ref{subsec:Uniform-bdd-condEHess}
and~\ref{subsec:Var-Hess} we can bound $a'$ and $b'$ independently
of $\alpha$ and $N$ provided $1\gg\alpha\gg N^{-1}$. The Lipschitz
constant $L$ can be bounded using Sections~\ref{subsec:Lipschitz-EHessf}
and~\ref{subsec:Lips-const-Hessf} to be $\mathcal{O}(1+1/(\sqrt{N\alpha}\,\alpha))$.
However after $\alpha\gets\beta\alpha$, $\widehat{f}_{\alpha}(\bm{w};\mathcal{D})-\inf_{\bm{w}'}\widehat{f}_{\alpha}(\bm{w}',\mathcal{D})$
is $\mathcal{O}(\alpha^{2})$: from Sections~\ref{subsec:Variance-of-gradients}
and~\ref{subsec:Changing-alpha}, $\nabla\widehat{f}_{\alpha}(\bm{w};\mathcal{D})=\mathcal{O}(\alpha)$.
Since there are bounds $0\prec a'I\preceq\Hess\widehat{f}_{\alpha}(\bm{w};\mathcal{D})\preceq b'I$
we can infer that $\widehat{f}_{\alpha}(\bm{w};\mathcal{D})-\inf_{\bm{w}'}\widehat{f}_{\alpha}(\bm{w}',\mathcal{D})$
is $\mathcal{O}(\alpha^{2})$. Thus $\gamma(\widehat{f}_{\alpha}(\bm{w};\mathcal{D})-\inf_{\bm{w}'}\widehat{f}_{\alpha}(\bm{w}',\mathcal{D}))=\mathcal{O}((1+1/(\sqrt{N\alpha}\,\alpha))\alpha^{2})=\mathcal{O}(\alpha^{2}+\sqrt{\alpha/N})=\mathcal{O}(1)$
provided $1\gg\alpha\gg N^{-1}$. 

Also, since $\epsilon=\eta\alpha$ in Algorithm~\ref{alg:SVM-l1-penalty-1-1},
line~6, $\epsilon_{0}/\epsilon=\mathcal{O}(1/(\alpha^{3}\eta))$
and so $\log\log(\epsilon_{0}/\epsilon)=\mathcal{O}(\log\log(1/\alpha))$.
This is very slowly growing function of $1/\alpha$. In fact, we can
show that only one Newton step is needed if $1\gg\alpha\gg N^{-1}$:
if $\left\Vert \nabla\widehat{f}_{\alpha}(\bm{w};\mathcal{D})\right\Vert \leq\frac{1}{2}(a')^{2}/L=\Omega(\sqrt{N\alpha}\,\alpha)$
then if $\bm{w}^{+}$ is the result of one full step of Newton's
method we have $(L/(2(a')^{2}))\left\Vert \nabla\widehat{f}_{\alpha}(\bm{w}^{+},\mathcal{D})\right\Vert \leq\left[(L/(2(a')^{2}))\left\Vert \nabla\widehat{f}_{\alpha}(\bm{w};\mathcal{D})\right\Vert \right]^{2}$.
This gives $\left\Vert \nabla\widehat{f}_{\alpha}(\bm{w}^{+},\mathcal{D})\right\Vert \leq L^{-1}\left[\mathcal{O}(L)\,\mathcal{O}(\alpha)\right]^{2}=\mathcal{O}(L\,\alpha^{2})=\mathcal{O}(\alpha/\sqrt{N\alpha})$.
If $\sqrt{N\alpha}\gg1/\eta$ then only one Newton step is needed. 
\end{proof}

\subsection{Endgame} \label{subsec:Endgame}
We now consider the asymptotic regime where $N\alpha\ll1$. This asymptotic
regime might not be needed in practice, as the values of $\alpha$
giving statistically reliable results will satisfy $N\alpha\gg1$.
However, for precisely identifying the support vectors, this is necessary.
In this regime, the line search may become necessary for Newton's
method to converge. The argument used here to show that only a bounded
number of Newton steps (with Armijo/backtracking line search) are
needed for each value of $\alpha$ is to essentially argue that the
functions are self-similar for $N\alpha\ll1$. Specifically, we identify
a limiting function $\bm{h}(\bm\omega)$ that is an approximation
to a scaled and shifted $\nabla\widehat{f}_{\alpha}(\bm{w})$.

Here we suppose that we have identified a set $\mathcal{I}$ of $m$
linearly independent data points $\bm{x}_{i}\in\R^{m}$ where
$(\bm{w}^{*})^{T}\bm{x}_{i}y_{i}=1$ for all $i\in\mathcal{I}$
where $\bm{w}^{*}$ is the minimizer of the unsmoothed SVM
objective function. We also assume that $\left|(\bm{w}^{*})^{T}\bm{x}_{j}y_{j}-1\right|\gg\alpha$
for all $j\not\in\mathcal{I}$. This would normally only happen when
$\alpha\ll1/N$ where $N$ is the number of data points.
\begin{lemma}
Suppose that $\bm{w}^{*}$ minimizes (\ref{eq:SVM-soft-margin})
and that the set $\mathcal{I}=\left\{ \,i\mid1=y_{i}\bm{x}_{i}^{T}\bm{w}^{*}\,\right\} $
has $m$ indexes, and that the vectors $\bm{x}_{i}$ for $i\in\mathcal{I}$
are linearly independent. Further assume that the unique solution
$\bm{\theta}$ of $\sum_{i\in\mathcal{I}}\theta_{i}y_{i}\bm{x}_{i}+\sum_{j\not\in\mathcal{I}\,\&\,1>y_{j}\bm{x}_{j}^{T}\bm{w}^{*}}y_{j}\bm{x}_{j}+N\lambda\bm{w}^{*}=\bm{0}$
has $0<\theta_{i}<1$ for all~$i\in\mathcal{I}$. If, in addition,
$1\gg\eta$ in the stopping criterion, the number of guarded Newton
steps needed for each value of $\alpha$ with $N^{-1}\gg\alpha$ is
asymptotically independent of $\alpha$, but dependent on $\beta$
and the index set $\mathcal{I}$.
\end{lemma}

\begin{proof}
(Outline) The optimality conditions for $\bm{w}=\bm{w}^{*}$
are
\begin{align*}
\bm{0} & \in N^{-1}\sum_{i\in\mathcal{I}}[0,1]\,y_{i}\bm{x}_{i}+\bm{g}^{*}\qquad\text{where}\\
\bm{g}^{*} & =N^{-1}\sum_{j\not\in\mathcal{I}}H(1-y_{j}(\bm{w}^{*})^{T}\bm{x}_{j})y_{j}\bm{x}_{j}+\lambda\bm{w}^{*},
\end{align*}
where $H(u)=1$ if $u>0$ and $H(u)=0$ otherwise. By linear independence
of the $\bm{x}_{i}$ over $i\in\mathcal{I}$, there are unique
$\theta_{i}\in[0,1]$ for $i\in\mathcal{I}$ where
\begin{equation}
-\bm{g}^{*}=N^{-1}\sum_{i\in\mathcal{I}}\theta_{i}y_{i}\bm{x}_{i}.\label{eq:g-star-theta}
\end{equation}
With linear independence of the $\bm{x}_{i}$'s, the $\theta_{i}$'s
are unique. For non-degeneracy, we assume that $\theta_{i}\neq0,\,1$
for all $i\in\mathcal{I}$. Then there are $\tau_{i}$ where $\psi'(\tau_{i})=y_{i}\theta_{i}$.
Let $\bm{d}$ be the unique solution of $\bm{d}^{T}\bm{x}_{i}=y_{i}\tau_{i}$
for $i\in\mathcal{I}$; there is a unique solution as $\left\{ \,\bm{x}_{i}\mid i\in\mathcal{I}\,\right\} $
is a basis for $\R^{m}$. Note that 
\[
\psi'\left(\frac{1-y_{j}\bm{w}^{T}\bm{x}_{j}}{\alpha}\right)=H(1-y_{j}\bm{w}^{T}\bm{x}_{j})+\mathcal{O}(\alpha^{2})\qquad\text{for }j\not\in\mathcal{I},
\]
since $\psi'(u)=H(u)+\mathcal{O}(u^{-2})$ as $u\to\pm\infty$, provided
$\left\Vert \bm{w}^{*}-\bm{w}\right\Vert =\mathcal{O}(\alpha)$.
Let 
\[
\bm{g}_{\alpha}(\bm{w})=-N^{-1}\sum_{j\not\in\mathcal{I}}\psi'\left(\frac{1-y_{j}\bm{w}^{T}\bm{x}_{j}}{\alpha}\right)y_{j}\bm{x}_{j}+\lambda\bm{w}.
\]
Now for $\left\Vert \bm{w}^{*}-\bm{w}\right\Vert =\mathcal{O}(\alpha)$,
$\bm{g}_{\alpha}(\bm{w})-\bm{g}^{*}=\mathcal{O}(\alpha^{2})+\lambda(\bm{w}-\bm{w}^{*})$.

If $\bm\omega=(\bm{w}^{*}-\bm{w})/\alpha$ then 
\begin{align*}
N\,\nabla\widehat{f}_{\alpha}(\bm{w};\mathcal{D}) & =-\sum_{i\in\mathcal{I}}\psi'(y_{i}\bm\omega^{T}\bm{x}_{i})y_{i}\bm{x}_{i}+N(\bm{g}_{\alpha}(\bm{w})-\bm{g}^{*})+N\bm{g}^{*}\\
 & =-\sum_{i\in\mathcal{I}}\psi'(y_{i}\bm\omega^{T}\bm{x}_{i})y_{i}\bm{x}_{i}+\sum_{i\in\mathcal{I}}\theta_{i}y_{i}\bm{x}_{i}+\mathcal{O}(N\alpha^{2})+N\alpha\,\lambda\bm\omega\\
 & \to-\sum_{i\in\mathcal{I}}\left[\psi'(y_{i}\bm\omega^{T}\bm{x}_{i})y_{i}-\theta_{i}\right]\bm{x}_{i}\qquad\text{as }N\alpha\downarrow0\text{ for fixed }\bm\omega.
\end{align*}
Since Newton's method with Armijo/backtracking line search is invariant
under translation and scaling, we expect the behavior of Newton's
method in the regime where $N\alpha\ll1$ to be consistent with Newton's
method applied to solving
\[
\bm{0}=\bm{h}(\bm\omega):=\sum_{i\in\mathcal{I}}\left[\psi'(y_{i}\bm{x}_{i}^{T}\bm{\omega})-\psi'(\tau_{i})\right]y_{i}\bm{x}_{i}.
\]
Note that $\bm{h}(\bm\omega)=\lim_{\alpha\downarrow0}-N\nabla\widehat{f}_{\alpha}(\bm{w}^{*}-\alpha\bm\omega)$.
Also note that $\bm{h}(\bm{\omega})=\nabla\widetilde{f}(\bm{\omega})$
with $\widetilde{f}(\bm{\omega})=\sum_{i\in\mathcal{I}}\left[\psi(y_{i}\bm{x}_{i}^{T}\bm{\omega})-\psi'(\tau_{i})y_{i}\bm{x}_{i}^{T}\bm{\omega}\right]$
making $\widetilde{f}$ a strictly convex and smooth function. Also,
$\widetilde{f}(\bm{\omega})-\widetilde{f}(\bm{d})=\lim_{\alpha\downarrow0}N\left(\widehat{f}_{\alpha}(\bm{w}^{*}-\alpha\bm{\omega})-\widehat{f}_{\alpha}(\bm{w}^{*})\right)$.

The solution of this limiting equation $\bm{h}(\bm\omega)=\bm{0}$
is clearly $\bm\omega^{*}=\bm{d}$. Let $\widehat{\bm{\omega}}\approx\bm{d}$
be the approximate solution from the previous guarded Newton iteration.
Then after reducing $\alpha$ by a factor of $\beta$, the new Newton
iteration for $\bm{\omega}$ has the starting point $\bm{\omega}_{0}=\widehat{\bm{\omega}}/\beta\approx\bm{d}/\beta$.
The smaller the value of $\eta$ the closer $\widehat{\bm{\omega}}$
is to $\bm{d}$ and the closer $\bm{\omega}_{0}$
is to $\bm{d}/\beta$. Consider the Jacobian matrix of $\bm{h}$:
\[
\nabla\bm{h}(\bm\omega)=\sum_{i\in\mathcal{I}}\psi''(y_{i}\bm{x}_{i}^{T}\bm\omega)\bm{x}_{i}\bm{x}_{i}^{T}=\Hess\widetilde{f}(\bm{\omega}).
\]
Essentially the guarded Newton iteration for $\bm{w}$ using
$\widehat{f}_{\alpha}(\bm{w};\mathcal{D})$ with $N^{-1}\gg\alpha>0$
approaches the guarded Newton iteration for $\bm{\omega}$
using $\widetilde{f}(\bm{\omega})$ with $\bm{\omega}_{k}\approx(\bm{w}^{*}-\bm{w}_{k})/\alpha$.
Thus the behavior of the iterates $\bm{w}_{k}$ can be represented
by the iterates $\bm{\omega}_{k}$ of the guarded Newton method
for $\widetilde{f}(\bm{\omega})$. For $\eta\ll1$ the starting
point is $\bm{\omega}_{0}\approx\bm{d}/\beta$. Applying
the guarded Newton method gives a sequence of iterations $\bm{\omega}_{k}$
for $k=0,1,2,\ldots$ where $\bm{\omega}_{k}\to\bm{\omega}^{*}$.
When $\nabla\widetilde{f}(\bm{\omega}_{k})^{T}\Hess\widetilde{f}(\bm{\omega}_{k})^{-1}\nabla\widetilde{f}(\bm{\omega}_{k})<\eta$
then the stopping condition is satisfied. The number of iterations
needed to achieve this stopping criterion for $\bm{\omega}_{k}$
is the asymptotic number of iterations needed between reductions in
$\alpha$ for $\bm{w}_{k}$. This number is independent of
$\alpha$. %
\end{proof}

If $1/\beta$ is large, then we could find that $\psi''(\tau_{i}/\beta)$
is small, and as a result, the Newton step is large compared to $\alpha$.
Different smoothing functions will change the asymptotics of this.
For example, if the smoothing function is $\psi(u)=\ln(1+e^{u})$,
then $\psi''(u)\sim e^{-\left|u\right|}$ so that $\left\Vert \nabla\bm{h}(\bm\omega)^{-1}\bm{h}(\bm\omega)\right\Vert $
grows exponentially in $\left\Vert \bm\omega\right\Vert $.%

\subsection{Summary of convergence results} \label{subsec:Summary-of-convergence}
In each of the asymptotic regimes (opening $\alpha\gg1$, midgame
$1\gg\alpha\gg N^{-1}$, and endgame $N^{-1}\gg\alpha$), we see that
only one step of the standard Newton method is typically sufficient
to replace $\alpha$ with $\beta\alpha$ for a fixed $\beta\in(0,1)$,
and no additional function evaluations are necessary. Between these
asymptotic regimes, convergence of the method is assured, although
the rate of convergence is not. While our analysis has not definitively
dealt with the transition cases, the authors believe that the method
can handle these cases without undue computational inefficiency. Thus
we expect that the number of passes over the data is $\mathcal{O}(\log(\alpha_{0}/\alpha_{\min}))$
where the method begins with $\alpha=\alpha_{0}$ and is terminated
when $\alpha<\alpha_{\min}$.

The use of an active set strategy may result in more passes over the
data in order to correctly identify the active set, but this strategy
also has the effect of keeping the number of active components of
$\bm{w}$ small. This means that the linear systems, and the
computed Hessian matrices, are both small. This avoids the communication
costs for transmitting large matrices across a cluster.

\section{Computational Results} \label{sec:results}
We use both real and synthetic data to compare SmSVM algorithms against well-known optimization algorithms such as the conjugate gradient and the stochastic gradient descent. We look at the ability of our models to accurately classify test data while maintaining, and in many cases improving, state-of-the-art training time. 

The rest of the section is organized as follows. Section \ref{subseq:data} describes the datasets used in the experiments. Section \ref{subseq:algo} presents the algorithms used for the performance comparison. In Section \ref{subseq:ncv}, we explain nested cross validation, which is used for hyperparameter tuning and for evaluating the performance of our algorithm. Finally, Section \ref{subseq:result} presents the results of the experiments.

\subsection{Data} \label{subseq:data}
The real datasets used in the experiments are called \textit{Australian Credit Approval}, \textit{Colon Cancer}, and \textit{Forest (Cover type)}, which are mostly acquired from the University of California Irvine Machine Learning Repository \cite{Dua2017}. A brief description of the data is given in Table \ref{table1}.

\begin{table}[H]
\tbl{Description of datasets used in performance comparison experiments.}
{\begin{tabular}{l r r r}
\toprule
\textbf{Name} & \textbf{Count} $(N)$ & \textbf{Dimension} $(m)$ & \textbf{Sparsity}\textsuperscript{a}  \\ \midrule
Australian Credit Approval & 690	& 14		& 20\%	 \\ 
Colon Cancer			& 62 		& 2,000	&  0\%	\\
Forest (Cover type) 		& 581,012 & 54		&  78\%	\\
Synthetic (tall)			& 10,000	& 50		& 0\%	\\
Synthetic (wide)		& 100	& 2,500	& 0\%	 \\ \bottomrule
\end{tabular}}
\tabnote{\textsuperscript{a}Sparsity refers to the percentage of the data with a value of $0$.}
\label{table1}
\end{table}

\subsubsection{Australian Credit Approval Data}
The famous Australian credit approval dataset, which originates from the StatLog project \cite{Mic94}, concerns credit card applications. The dataset consists of 14 feature variables and a class label that quantifies the approval decision. There are 690 instances in the dataset. Out of the 14 features, 6 are numerical and 8 are categorical. Note that four categorical variables have more than two categories. This means one-hot encoding is required for these features if they are nominal. Since all feature names and values have been changed for the sake of confidentiality of the data, we do not know if the categorical variables are nominal or ordinal. However, we observe a strong correlation between those variables and the approval decision. It leads us to assume that those categorical variables are ordinal. Experimental results are also in line with our assumption. The test accuracies without one-hot encoding are higher than those with one-hot encoding. Finally, it is worth noting that the dataset is chosen for its shape as it is tall with $N = 690$ instances and $m = 14$ features. The dataset is available on the UCI Machine Learning Repository website \cite{Dua2017}.

\subsubsection{Colon Cancer Data}
In this experiment, we use the colon cancer dataset \citep{Alo99} to classify tissues based on gene expression. The dataset\footnote{The dataset is available on Uri Alon's website at the following address with the name ``Affymetrix array data as described in the paper": \href{https://www.weizmann.ac.il/mcb/UriAlon/download/downloadable-data}{https://www.weizmann.ac.il/mcb/UriAlon/download/downloadable-data}} consists of 2000 feature variables and 62 gene expression. There are 22 normal and 40 tumor colon tissue samples in the dataset. As in many gene expression data, a number of features (genes) are highly correlated in the colon cancer dataset. Using correlation analysis by eliminating features with correlation greater than 0.7, we obtain a reduced set of data with 215 features. Next, standardization is performed for obvious reasons. We further use principal component analysis for dimensionality reduction. After the feature selection by correlation analysis and principal component analysis, only $20$ features remain. We finally apply our algorithm to the remaining data with $N=62$ instances and $m=20$ features.

\subsubsection{Forest (Cover Type) Data}
The forest cover-type dataset \citep{Bla98} is used to predict forest cover type from cartographic variables. The dataset comprises 54 features and 581,012 observations. Some features are elevation, aspect, slope or hillshade in the morning and in the afternoon. Note that 44 features are qualitative (all binary) and 10 are quantitative. The dataset has originally 7 classes (forest types). We modify it into a binary classification problem where the aim is to separate class 2 (the most common class) from the other 6 classes. The dataset is chosen due to its size: It has the largest number of observations among the datasets we run in our experiments. 

\subsubsection{Synthetic Data}
The synthetic data is generated by creating two centroids $\bm{c}_1, \bm{c}_2 \in \mathbb{R}^m$ with components randomly sampled from $\mathcal{N}(0,1)$, scaling the centroids, and then sampling instances from $\mathcal{N}(\bm{c}_i, I)$ for $i=1,2$. Each class has the same number of instances. We create two synthetic datasets with different shapes: 
\begin{itemize}
\item[(i)] synthetic tall dataset with $m = 50$ and $N = 10,000$;
\item[(ii)] synthetic wide dataset with $m = 2500$ and $N = 100$.
\end{itemize}

\subsection{Algorithms} \label{subseq:algo}
We compare our algorithms against the conjugate gradient, stochastic gradient descent, and coordinate descent (via LIBLINEAR \cite{Fan2008}). Table \ref{table2} summarizes the algorithms and their naming convention. All algorithms with $\ell^2$ regularization solve the optimization problem defined by equation \eqref{eq:SVM-soft-margin}, while those with $\ell^1$ and $\ell^2$ regularization minimize the loss function defined in equation \eqref{eq:L1-SVM}. In the case of conjugate gradient and stochastic gradient descent, since our loss function is non-smooth, we choose an arbitrary subgradient in place of the gradient, where the subgradient is defined as any element of the subdifferential,
\[
\partial f(\bm{x}) = \left\{ \bm{g} \in \mathbb{R}^n \mid f(\bm{y}) \geq f(\bm{x}) + \bm{g}^T(\bm{y}-\bm{x}), \quad \forall \bm{y} \in \R^n \right\}.
\]
Since the hinge-loss function and $\ell^1$-norm are convex, the behavior is sufficiently close to that of conjugate gradient and stochastic gradient descent on a smooth function.

All the experiments for each algorithm are carried out in Julia v1.7.3. In the SvSVM algorithms, the value of $\beta$ used is $0.1$; that is the factor by which the smoothing parameter $\alpha$ is reduced when needed. We use $\alpha_{\min} = 10^{-5}$ for validation runs and $\alpha_{\min} = 10^{-6}$ for test runs.

As for the conjugate gradient, we implement Hager and Zhang's algorithm ``CG\_DESCENT, a conjugate gradient method with guaranteed descent" \cite{Hag2006} incorporated with \cite{Hag2013}. In this algorithm, the termination condition employs a ``unit-correct" expression in place of a condition on gradient components. Note that the line search used here is exactly the one proposed in Hager and Zhang \cite{Hag2006}. The experiments for the CG are all implemented in Julia using Optim.jl package \cite{Mog2018}.

The LinearSVC model, which is a Julia wrapper over LIBLINEAR \cite{Fan2008}, solves a scaled version of equation \eqref{eq:SVM-soft-margin}; that is,
\[
\frac{1}{2} \bm{w}^T \bm{w} + C \sum_{i=1}^N \max\left\{ 0,\ 1-y_i (\bm{w}^T \bm{x}_i + b) \right\}, \quad C>0.
\]
In our case, this is optimized via coordinate descent (see \cite{Fan2008} for details). Hyper-parameter tuning for $C$ is performed by nested cross-validation, which is explained in detail in Section \ref{subseq:ncv}. The LIBLINEAR implementation was accessed via LIBSVM.jl \cite{Pas}, which provides a Julia wrapper on the C++ coded SVM package LIBSVM \cite{Cha2011}.

\begin{table}
\tbl{Summary of algorithms and their naming conventions.}
{\begin{tabular}{ll}
\toprule
\textbf{Name} & \textbf{Description} \\ \midrule
SmSVM-$\ell^2$ 		& Our algorithm with $\ell^2$ regularization \\ 
SmSVM-$\ell^1$-$\ell^2$ 	& Our algorithm with $\ell^1$ and $\ell^2$ regularization \\ 
CG-$\ell^2$			& Conjugate gradient with $\ell^2$ regularization \\
CG-$\ell^1$-$\ell^2$		& Conjugate gradient with $\ell^1$ and $\ell^2$ regularization \\
LinearSVC			& LIBLINEAR \\
Nesterov-$\ell^2$		& Nesterov (SGD) optimizer with $\ell^2$ regularization \\ 
Nesterov-$\ell^1$-$\ell^2$	& Nesterov (SGD) optimizer with $\ell^1$ and $\ell^2$ regularization \\
AdaDelta-$\ell^2$ 		& AdaDelta (SGD) optimizer with $\ell^2$ regularization \\ 
AdaDelta-$\ell^1$-$\ell^2$	& AdaDelta (SGD) optimizer with $\ell^1$ and $\ell^2$ regularization \\ 
Adam-$\ell^2$ 			& Adam (SGD) optimizer with $\ell^2$ regularization \\ 
Adam-$\ell^1$-$\ell^2$	& Adam (SGD) optimizer with $\ell^1$ and $\ell^2$ regularization \\ \bottomrule
\end{tabular}}
\label{table2}
\end{table}

We consider three different versions of SGD in our experiments: Nesterov momentum \cite{Nes83}, AdaDelta \cite{Zei2012}, and Adam \cite{Ada2014}. The idea of momentum is incorporated to speed progress since the gradient descent takes a long time to traverse a nearly flat surface. Nesterov momentum \cite{Nes83} uses the gradients of projected positions in order to avoid possible overshoot of the minima at the bottom of a valley due to the acceleration of momentum. The values of parameters for Nesterov momentum used in the experiments are the learning rate $\eta = 0.001$ and the momentum decay $\rho=0.9$.

AdaDelta \cite{Zei2012} is a per-dimension learning rate method for SGD. The method simply adapts learning rates based on a decaying moving average of past squared gradients. The method is an extension of the AdaGrad \cite{Duc2011} and RMSProp \cite{Hin} algorithms that addresses the drawbacks of these methods. AdaDelta overcomes AdaGrad's primary weakness, the reduction in the effective learning rate during training. Furthermore, AdaDelta requires no manual tuning of a learning rate unlike RMSProp. The values of parameters used for AdaDelta are the gradient decay $\rho=0.9$ and the small constant $\epsilon = 10^{-8}$.  

Adam optimizer \cite{Ada2014} is another version of the SGD method that takes the advantage of both momentum and adaptive gradient methods. Adam is based on adaptive estimation of both first and second-order moments. According to Kingma and Ba \cite{Ada2014}, Adam is ``computationally efficient, has little memory requirement, invariant to diagonal rescaling of gradients, and is well suited for problems that are large in terms of data/parameters". The values of parameters used for Adam are the learning rate $\eta = 0.001$ along with the decay of momentums: $\beta_1 = 0.9$ (the first moment estimate) and $\beta_2 = 0.8$ (the second moment estimate). The Nesterov momentum, AdaDelta, and Adam optimizers are all implemented in Julia using Flux.jl package \cite{Inn2018}.

\subsection{Nested Cross Validation} \label{subseq:ncv}
Nested cross-validation is used simultaneously for hyperparameter tuning and for evaluating the performance of our algorithm. Cawley and Talbot \cite{Caw2010} demonstrates that over-fitting in model selection can lead to biased performance evaluation. In order to avoid the bias in performance evaluation, model selection must be viewed as an integral part of the model fitting procedure as in nested cross-validation. 

Nested cross-validation has a double loop: an outer loop serves for performance evaluation, and an inner loop serves for hyper-parameter selection. It consists of an outer loop of $k$ sets and an inner loop of $\ell$ sets. The dataset is divided into $k$ sets. In each iteration of outer loop, a set is chosen as the outer test set and the rest is combined into the corresponding outer training set. Each outer training set is further split into $\ell$ sets. Similarly, in each iteration of inner loop, a set is chosen as the validation set and the rest is combined into the corresponding inner training set. In the inner loop, best hyperparameters are selected to minimize error on the validation set. Using the hyperparameters selected in the inner loop, we train the weights and then test our model's prediction on the test dataset. In our estimation, we use $k = 10$ and $\ell = 6$ for all datasets except the forest (cover type) dataset. The values $k=5$ and $\ell=4$ are used for forest dataset due to computational requirements. 

Before each training, we standardize training and test/validation datasets. We first perform standardization over the training data, and then the validation/set dataset is standardized using the mean and standard deviation of training dataset. The algorithm used for nested cross-validation is expressed in pseudocode in Algorithm \ref{algo_nestedcv}.

\begin{algorithm}
\begin{algorithmic}[1]
\Require{$\mathcal{D} = (X,\bm{y})$, hyper-parameter space $\lambda \in \Lambda$, $\mu \in M$, number of folds $k, \ell$}
\State$\mathcal{D}_1, \ldots, \mathcal{D}_k \gets$ split $\mathcal{D}$ into $k$ sets
\For{$i = 1,\ldots,k$} \Comment{Outer loop for performance evaluation}
	\State$\mathcal{D}_{\text{test}} \gets \mathcal{D}_i$ \Comment{Test dataset}
	\State$\mathcal{D}_{\text{train}} \gets \mathcal{D} \setminus \mathcal{D}_i$ \Comment{Outer training dataset}
	\State$\mathcal{D}_{\text{train}}^1, \ldots, \mathcal{D}_{\text{train}}^\ell \gets$ split $\mathcal{D}_{\text{train}}$ into $\ell$ sets
	\For{$j = 1,\ldots,\ell$}\Comment{Inner loop for hyper-parameter tuning}
		\State$\mathcal{D}_{\text{val}} \gets \mathcal{D}_{\text{train}}^j$ \Comment{Validation dataset}
		\State$\mathcal{D}_{\text{inner}} \gets \mathcal{D}_{\text{train}} \setminus \mathcal{D}_{\text{train}}^j$ \Comment{Inner training dataset}
		\State Standardize $X_{\text{inner}}$ and $X_{\text{val}}$
		\For{$(\lambda, \mu) \in \Lambda \times M$}
			\State$\bm{w} \gets \text{\sc{SVMsmooth}}(X_{\text{inner}}, \bm{y}_{\text{inner}}, \lambda, \mu)$\;
			\State$\text{acc}_{\lambda,\mu}^j \gets$ calculate validation accuracy using $\bm{w}$ over $\mathcal{D}_{\text{val}}$
		\EndFor
	\EndFor
	\State$\text{acc}_{\lambda,\mu} \gets \frac{1}{\ell}\left(\text{acc}_{\lambda,\mu}^1 + \cdots + \text{acc}_{\lambda,\mu}^\ell \right)$ for all $(\lambda, \mu)$
	\State$(\lambda^*, \mu^*) \gets $ pick $(\lambda,\mu )$ with the highest validation accuracy, $\text{acc}_{\lambda,\mu}$
	\State Standardize $X_{\text{train}}$ and $X_{\text{test}}$
	\State $\bm{w} \gets \text{\sc{SVMsmooth}}(X_{\text{train}}, \bm{y}_{\text{train}}, \lambda^*, \mu^*)$
	\State $\text{acc}_i \gets$ calculate test accuracy using $\bm{w}$ over $\mathcal{D}_{\text{test}}$
\EndFor
\State $\mathbf{return}\; \frac{1}{k} \left( \text{acc}_1 + \cdots + \text{acc}_k \right)$\;
\end{algorithmic}
\caption{Nested cross-validation} \label{algo_nestedcv} 
\end{algorithm}

\subsection{Results} \label{subseq:result}
As seen in Table~\ref{tab:results}, SmSVM\textendash $\ell^{1}$\textendash $\ell^{2}$ performs well across a variety of dataset types. This may be due to the feature selection property of the $\ell^{1}$ norm. We optimize the matrix-vector and vector-vector operations by reducing the problem size to that of the active set dimension. The reduction in problem size yields substantial computational savings in problems where the active-set is small. In the synthetic wide dataset, our algorithms shows a relatively poor performance in training time. Thus, the SmSVM algorithms may take some time to finish training when the number of features is much greater than the number of samples, i.e. $m \gg N$.


Perhaps the most surprising result is the performance on the forest (cover type) dataset. Consisting of nearly 600,000 data points which takes up roughly 140MB, LinearSVC \cite{Fan2008} took nearly a minute to train on this dataset, achieving a best-in-class test accuracy, while SmSVM\textendash $\ell^{2}$ trained in just about twenty seconds and achieving second-place test accuracy. The fastest algorithms in terms of training time for the forest dataset are the SGD variants.

The results of the numerical tests do not have clear winners and losers, except perhaps that the conjugate gradient method performs worse than the other methods. For the most of the datasets, one of the SmSVM versions (with either the $\ell^{2}$ or both $\ell^{2}$ and $\ell^{1}$ penalties) is better than the other methods on \emph{either} the computational time or the test accuracy. We therefore claim that SmSVM is a competitive method, although there are modifications to SmSVM that have the potential to improve its performance. The point to be made is that second order methods like Newton's method should not be disregarded even for large datasets. 

\renewcommand{\arraystretch}{1.2}
\begin{table}
\tbl{Numerical results for all datasets.}
{
\begin{tabular}{l  cc  cc  cr  cc  cr}
\toprule
\multirow{2}{*}{\textbf{Algorithm}} & \multicolumn{2}{c}{\textbf{Australian}} & \multicolumn{2}{c}{\textbf{Colon Cancer}} &  \multicolumn{2}{c}{\textbf{Covertype}} &  \multicolumn{2}{c}{\textbf{Synthetic (tall)}} &  \multicolumn{2}{c}{\textbf{Synthetic (wide)}}  \\ 
\cmidrule(lr){2-3}
\cmidrule(lr){4-5}
\cmidrule(lr){6-7}
\cmidrule(lr){8-9}
\cmidrule(lr){10-11}
& Acc & Time (s) & Acc & Time (s) & Acc & Time (s) & Acc & Time (s) & Acc & Time (s) \\ \midrule
SmSVM-$\ell^2$ 		& 86.67 & 0.0128 & 79.03 & 0.00212 & 76.32 & 22.52 & 100.0 & 0.2272 & 100.0 & 13.4500  \\ 
SmSVM-$\ell^1$-$\ell^2$	& 86.67 & 0.0124 & 80.65 & 0.00174 & 76.32 & 35.59 & 100.0 & 0.1981 & 100.0 & 7.6830 \\
CG-$\ell^2$ 			& 85.51 & 0.1005 & 74.19 & 0.03580 & 76.32 & 193.35 & 100.0 & 0.4886 & 100.0 & 0.0740 \\
CG-$\ell^1$-$\ell^2$ 		& 86.52 & 0.0672 & 77.42 & 0.02972 & 76.32 & 146.23 & 100.0 & 0.1779 & 100.0 & 0.0661 \\
LinearSVC 			& 86.81 & 0.0002 & 74.19 & 0.00008 & 76.35 &	 51.57 & 100.0 & 0.0071 & 100.0 & 0.0028 \\
Nesterov-$\ell^2$ 		& 86.67 & 0.0693 & 77.42 & 0.00146 & 74.57 & 14.85 & 100.0 & 0.2581 & 100.0 & 0.0169 \\
Nesterov-$\ell^1$-$\ell^2$ & 86.38 & 0.0810 & 79.03 & 0.00194 & 74.52 & 19.38 & 100.0 & 0.3622 & 100.0 & 0.0215 \\ 
AdaDelta-$\ell^2$ 		& 86.96 & 0.0170 & 79.03 & 0.00144 & 76.29 & 14.97 & 100.0 & 0.2588 & 100.0 & 0.0175 \\
AdaDelta-$\ell^1$-$\ell^2$ & 87.39 & 0.0202 & 80.65 & 0.00180 & 76.29 & 19.05 & 100.0 & 0.3450 & 100.0 & 0.0219 \\
Adam-$\ell^2$ 			& 85.65 & 0.0686 &79.03 & 0.00145 & 75.76 & 14.44 & 100.0 & 0.2477 & 100.0 & 0.0182  \\
Adam-$\ell^1$-$\ell^2$ 	& 86.38 & 0.0207 & 80.65 & 0.00169 & 75.86 & 18.76 & 100.0 & 0.3526 & 100.0 & 0.0217  \\
\bottomrule
\end{tabular}
}
\label{tab:results}
\end{table}

\section{Conclusion} \label{sec:Conclusion}
We have introduced SmSVM, a new approach to solving soft-margin SVM, which is capable of strong test accuracy without sacrificing training speed. This is achieved by smoothing the hinge-loss function and using an active set approach to the the $\ell^{1}$ penalty. SmSVM provides improved test accuracy over the other methods with comparable, and in some cases reduced, training time. SmSVM uses many fewer gradient calculations and a modest number of passes over the data to achieve its results, meaning it scales well for large datasets. SmSVM\textendash $\ell^{1}$\textendash $\ell^{2}$ optimizes its matrix-vector and vector-vector calculations by reducing the problem size to that of the active set. For even modestly sized problems, this results in significant savings with respect to computational complexity. 

Overall the results are quite promising. On both real and synthetic datasets, our algorithms outperform or tie the competition in test accuracy most of the time. Moreover, for tall datasets (number of samples $\gg$ number of features), SmSVM algorithms are among the best in terms of training time. The time savings are increasingly significant as the number of data points grows. 

The convergence analysis leads to a number of questions: How does the number of passes over the data change we make $\beta$ small? Is there a near optimal value of $\beta$ that is independent of the dataset? While the method is designed to operate more efficiently with large datasets ($N$ large), we can control the value of $N$ by taking statistically representative subsets of the dataset; we could select a subset of $N_{1}\ll N$ data points from the original dataset, then solve the smoothed problem taking $\alpha$ to approximately $m/N_{1}$, then restarting the algorithm on a larger subset of $N_{2}\gg N_{1}$ data points with this value of $\alpha$, then reducing $\alpha$ down to approximately $m/N_{2}$, etc. Rapid convergence would still be expected, even though the computations are performed on a smaller dataset until the last steps using the full dataset. In this way, the computational cost can be reduced to the equivalent of just a few passes over the full dataset.

SmSVM is implemented in Julia, making it easy to modify and understand. The use of Julia keeps linear algebra operations optimized. This is important when competing against frameworks such as LIBLINEAR, which is implemented in C++. Testing SmSVM on larger datasets, incorporating GPU acceleration to the linear algebra, and exploring distributed implementations are promising future directions.


\begin{thebibliography}{}

\end{thebibliography}


\begin{thebibliography}{99}

\bibitem{Ada2014}
D.P. Kingma and J. Ba, \emph{Adam: A method for stochastic optimization}, arXiv preprint (2014). Available at \href{https://arxiv.org/pdf/1412.6980}{https://arxiv.org/pdf/1412.6980}.

\bibitem{Alo99}
U. Alon, N. Barkai, D.A. Notterman, K. Gish, S. Ybarra, D. Mack, and A.J. Levine, \emph{Broad patterns of gene expression revealed by clustering analysis of tumor and normal colon tissues probed by oligonucleotide arrays}, Proceedings of the National Academy of Sciences 96.12 (1999), pp. 6745--6750.

\bibitem{Arm66}
L. Armijo, \emph{Minimization of functions having Lipschitz continuous first partial derivatives}, Pacific J. Math. 16 (1966), pp. 1--3.

\bibitem{Bla98}
J.A. Blackard, \emph{Comparison of neural networks and discriminant analysis in predicting forest cover types}, Ph.D. diss., Colorado State University, 1998.

\bibitem{Ble96}
G.E. Blelloch, \emph{Programming parallel algorithms}, Communications of the ACM 39.3 (1996), pp. 85--97.

\bibitem{Boy2004}
S. Boyd, S.P. Boyd, and L. Vandenberghe, \emph{Convex Optimization}, Cambridge University Press, Cambridge, 2004.

\bibitem{Caw2010}
G.C. Cawley and N.L. Talbot, \emph{On over-fitting in model selection and subsequent selection bias in performance evaluation}, J. Mach. Learn. Res. 11 (2010), pp. 2079--2107.

\bibitem{Cha2011}
C.C. Chang and C.J. Lin, \emph{LIBSVM: A library for support vector machines}, ACM Transactions on Intelligent Systems and Technology 2.3 (2011), pp. 1--27.

\bibitem{Dua2017}
D. Dua and C. Graff, UCI Machine Learning Repository, 2017; dataset available at \href{http://archive.ics.uci.edu/ml}{http://archive.ics.uci.edu/ml}.

\bibitem{Duc2011}
J. Duchi, E. Hazan, and Y. Singer, \emph{Adaptive subgradient methods for online learning and stochastic optimization}, J. Mach. Learn. Res. 12.7 (2011).

\bibitem{Dud2000}
R.O. Duda, P.E. Hart, and D.G. Stork, \emph{Pattern Classification}, 2nd ed., Wiley-Interscience, 2000.

\bibitem{Fan2008}
R.E. Fan, K.W. Chang, C.J. Hsieh, X.R. Wang, and C.J. Lin, \emph{LIBLINEAR: A library for large linear classification}, J. Mach. Learn. Res. 9 (2008), pp. 1871--1874.

\bibitem{Fle87}
R. Fletcher, \emph{Practical Methods of Optimization}, 2nd ed., A Wiley-Interscience Publication, John Wiley \& Sons Ltd., Chichester, 1987.

\bibitem{Hag2006}
W.W. Hager and H. Zhang, \emph{Algorithm 851: CG\_DESCENT, a conjugate gradient method with guaranteed descent}, ACM Trans. Math. Software 32.1, (2006), pp. 113--137.

\bibitem{Hag2013}
W.W. Hager and H. Zhang, \emph{The limited memory conjugate gradient method}, SIAM Journal on Optimization 23.4 (2013), pp. 2150--2168.

\bibitem{Haj2018}
J. Hajewski, S. Oliveira, and D. Stewart, \emph{Smoothed Hinge Loss and $\ell^1$ Support Vector Machines}, IEEE International Conference on Data Mining Workshops (ICDMW), 2018, pp. 1217--1223.

\bibitem{Hin}
G. Hinton, \emph{RMSProp: Divide the gradient by a running average of its recent magnitude}, Coursera Lecture 6e. Available at \href{http://www.cs.toronto.edu/~tijmen/csc321/slides/lecture_slides_lec6.pdf}{http://www.cs.toronto.edu/~tijmen/csc321/slides/lecture\_slides\_lec6.pdf}.

\bibitem{Hor85}
R.A. Horne and C.A. Johnson, \emph{Matrix Analysis}, Cambridge Uni. Press, Cambridge, 1985.

\bibitem{Inn2018}
M. Innes, \emph{Flux: Elegant Machine Learning with Julia}, Journal of Open Source Software 3.25 (2018).

\bibitem{Kan49}
L.V. Kantorovi\u{c}, \emph{On Newton's method}, Trudy Mat. Inst. Steklov. 28 (1949), pp. 104--144.

\bibitem{Lie2001}
E.H. Lieb and M. Loss, \emph{Analysis}, 2nd ed., Graduate Studies in Mathematics Vol. 14, American Mathematical Society, Providence, RI, 2001.

\bibitem{Loe77}
M. Lo\`{e}ve, \emph{Probability Theory I}, 4th ed., Springer-Verlag, New York, NY, 1977.

\bibitem{Mic94}
D. Michie, D. Spiegelhalter, and C. Toylor, \emph{Machine Learning, Neural Network and Statistical Classification}, Ellis Horwood Series in Artificial Intelligence, 1994.

\bibitem{Mog2018}
P.K. Mogensen and A.N. Riseth, \emph{Optim: A mathematical optimization package for Julia}, Journal of Open Source Software 3.24 (2018).

\bibitem{Nes83}
Y.E. Nesterov, \emph{A method for solving the convex programming problem with convergence rate $O\bigl(\frac{1}{k^ 2} \bigr)$}, In Dokl. Akad. Nauk SSSR, Vol. 269, (1983), pp. 543--547.

\bibitem{Nes94}
Y. Nesterov and A. Nemirovskii,\emph{Interior-point polynomial algorithms in convex programming}, SIAM Studies in Applied Mathematics Vol. 13, Society for Industrial and Applied Mathematics (SIAM), Philadelphia, PA, 1994.

\bibitem{Niu2011}
F. Niu, B. Recht, C. Re, and S.J. Wright, \emph{Hogwild!: A Lock-free Approach to Parallelizing Stochastic Gradient Descent}, Advances in Neural Information Processing Systems 24, 2011.

\bibitem{Noc2006}
J. Nocedal and S.J. Wright, \emph{Numerical Optimization}, 2nd ed., Springer, NewYork, 2006.

\bibitem{Osb2000}
M.R. Osborne, B. Presnell, and B.A. Turlach, \emph{A new approach to variable selection in least squares problems}, IMA J. Numer. Anal. 20 (2000), pp. 389--403. 

\bibitem{Pas}
M. Pastell, \emph{The LIBSVM.jl library}; software available at \href{https://github.com/mpastell/LIBSVM.jl}{https://github.com/mpastell/LIBSVM.jl}.

\bibitem{Rob51}
H. Robbins and S. Monro, \emph{A stochastic approximation method}, The Annals of Mathematical Statistics (1951), pp. 400--407.

\bibitem{Zei2012}
M.D. Zeiler, \emph{Adadelta: an adaptive learning rate method}, arXiv preprint (2012). Available at \href{https://arxiv.org/abs/1212.5701}{https://arxiv.org/abs/1212.5701}.

\bibitem{Zhu2003}
J. Zhu, S. Rosset, T. Hastie, and R. Tibshirani, \emph{1-norm Support Vector Machines}, Advances in Neural Information Processing Systems 16, 2003.

\end{thebibliography}
\end{document}